# A Linearly Convergent Majorized ADMM with Indefinite Proximal Terms for Convex Composite Programming and Its Applications

Ning Zhang,* Jia Wu,† Liwei Zhang ‡

February 7, 2018

**Abstract** This paper aims to study a majorized alternating direction method of multipliers with indefinite proximal terms (iPADMM) for convex composite optimization problems. We show that the majorized iPADMM for 2-block convex optimization problems converges globally under weaker conditions than those used in the literature and exhibits a linear convergence rate under a local error bound condition. Based on these, we establish the linear rate convergence results for a symmetric Gaussian-Seidel based majorized iPADMM, which is designed for multi-block composite convex optimization problems. Moreover, we apply the majorized iPADMM to solve different types of regularized logistic regression problems. The numerical results on both synthetic and real datasets demonstrate the efficiency of the majorized iPADMM and also illustrate the effectiveness of the introduced indefinite proximal terms.

## 1 Introduction

In this paper, we consider the following convex composite optimization problem:

$$\begin{aligned}\min_{y,z} \quad & p(y) + f(y) + q(z) + g(z) \\ \text{s.t.} \quad & \mathcal{A}^* y + \mathcal{B}^* z = c, \\ & y \in \mathcal{Y}, z \in \mathcal{Z},\end{aligned} \qquad (1.1)$$

where $\mathcal{X}, \mathcal{Y}$ and $\mathcal{Z}$ are given finite dimensional Euclidean spaces each equipped with an inner product $\langle \cdot, \cdot \rangle$ and its induced norm $\|\cdot\|$; $f : \mathcal{Y} \to (-\infty, +\infty)$ and $g : \mathcal{Z} \to (-\infty, +\infty)$ are two convex functions with Lipschitz continuous gradients on $\mathcal{Y}$ and $\mathcal{Z}$, respectively; $p : \mathcal{Y} \to (-\infty, +\infty]$ and $q : \mathcal{Z} \to (-\infty, +\infty]$ are two closed proper convex (not necessarily smooth) functions; $\mathcal{A}^* : \mathcal{Y} \to \mathcal{X}$ and $\mathcal{B}^* : \mathcal{Z} \to \mathcal{X}$ are adjoints of the linear operators $\mathcal{A} : \mathcal{X} \to \mathcal{Y}$ and $\mathcal{B} : \mathcal{X} \to \mathcal{Z}$, respectively; and $c \in \mathcal{X}$.

---
*Future Resilient Systems, Singapore-ETH center; Department of Mathematics, National University of Singapore, Singapore. (matzhn@nus.edu.sg). The Singapore-ETH Centre (SEC) was established as a collaboration between ETH Zurich and National Research Foundation (NRF) Singapore (FI 370074011) under the auspices of the NRF's Campus for Research Excellence and Technological Enterprise (CREATE) programme.

†School of Mathematical Sciences, Dalian University of Technology, China. (wujia@dlut.edu.cn).

‡School of Mathematical Sciences, Dalian University of Technology, China. (lwzhang@dlut.edu.cn). This author's research is supported by the National Natural Science Foundation of China under project 11571059 and 11731013.

Let $\mathcal{M} : \mathcal{X} \to \mathcal{X}$ be a self-adjoint linear operator (not necessarily semidefinite positive), denote $\|x\|_\mathcal{M}^2 := \langle x, \mathcal{M}x \rangle$. Let $\sigma > 0$ be a given parameter. The augmented Lagrangian function of (1.1) is defined by

$$\mathcal{L}_\sigma(y, z; x) = f(y) + p(y) + g(z) + q(z) + \langle x, \mathcal{A}^* x + \mathcal{B}^* y - c \rangle + \frac{\sigma}{2} \|\mathcal{A}^* x + \mathcal{B}^* y - c\|^2.$$

Consider the following general 2-block ADMM iterative scheme,

$$\begin{cases} y^{k+1} \in \mathrm{argmin} \left\{ \mathcal{L}_\sigma(y, z^k; x^k) + \tfrac{1}{2} \|y - y^k\|_\mathcal{S}^2 \mid y \in \mathcal{Y} \right\}, \\ z^{k+1} \in \mathrm{argmin} \left\{ \mathcal{L}_\sigma(y^{k+1}, z; x^k) + \tfrac{1}{2} \|z - z^k\|_\mathcal{T}^2 \mid z \in \mathcal{Z} \right\}, \\ x^{k+1} = x^k + \tau\sigma(\mathcal{A}^* y^{k+1} + \mathcal{B}^* z^{k+1} - c). \end{cases} \quad (1.2)$$

It is well known that if $\mathcal{S} = 0$ and $\mathcal{T} = 0$, the iterative scheme (1.2) is exactly the classic ADMM designed by Glowinski and Marroco [22] and Gabay and Mercier [21]; if $\mathcal{S} \succ 0$, $\mathcal{T} \succ 0$ and $\tau = 1$, iterative scheme (1.2) reduces to the method of the proximal ADMM introduced by Eckstein [14]; if both $\mathcal{S}$ and $\mathcal{T}$ are self-adjoint positive semidefinite linear operators, $\tau \in (0, (1 + \sqrt{5})/2)$, iterative scheme (1.2) is known as the semiproximal ADMM (sPADMM) which is proposed by Fazel et al. [17]. To know more about the above mentioned works and their relationships with well known methods, such as proximal point algorithm (PPA) and Douglas-Rachford (DR) splitting method, we refer the readers to [22, 18, 16, 15, 14, 23, 7, 24, 34, 27].

One of the most important motivations behind ADMM is to fully use the separable structures in the problems. In other words, a potential assumption of using ADMM is that each subproblem can be efficiently solved. Generally speaking, if $f$ or $g$ is not a quadratic or linear function, the corresponding subproblem does not have closed-form solutions or cannot be solved easily. In order to continue enjoying benefits of the separable structure, Li et al. [34] extended the sPADMM to a majorized ADMM with indefinite proximal terms (iPADMM). Compared with the majorized techniques mentioned in [7, 1], the majorized iPADMM uses the positive semidefinite operators $\widehat{\Sigma}_f$ and $\widehat{\Sigma}_g$ (see (1.5) and (1.6)) instead of the Lipschitz constants of the gradient mappings $\nabla f$ and $\nabla g$. The motivation behind this is for the better numerical performance. This will be illustrated in Table 2. Li et al. [34] established the global convergence and the iteration-complexity in the nonergodic sense of the majorized iPADMM, but not the rate of convergence.

In this paper, we further study the convergence and the rate of convergence of the majorized iPADMM presented in [34]. Now, we recall the majorized iPADMM. Since both $f(\cdot)$ and $g(\cdot)$ are smooth convex functions, there exist self-adjoint and positive semidefinite linear operators $\Sigma_f$ and $\Sigma_g$ such that for any $y, y' \in \mathcal{Y}$ and any $z, z' \in \mathcal{Z}$,

$$f(y) \geq f(y') + \langle \nabla f(y'), y - y' \rangle + \frac{1}{2} \|y - y'\|_{\Sigma_f}^2, \quad (1.3)$$

$$g(z) \geq g(z') + \langle \nabla g(z'), z - z' \rangle + \frac{1}{2} \|z - z'\|_{\Sigma_g}^2. \quad (1.4)$$

In addition, by the condition that the gradients $\nabla f(\cdot)$ and $\nabla g(\cdot)$ are Lipschitz continuous, we know that there exist self-adjoint and positive semidefinite linear operators $\widehat{\Sigma}_f \succeq \Sigma_f$ and $\widehat{\Sigma}_g \succeq \Sigma_g$ such that for any $y, y' \in \mathcal{Y}$ and any $z, z' \in \mathcal{Z}$,

$$f(y) \le \hat{f}(y, y') := f(y') + \langle \nabla f(y'), y - y' \rangle + \frac{1}{2}\|y - y'\|^2_{\widehat{\Sigma}_f}, \tag{1.5}$$

$$g(z) \le \hat{g}(z, z') := g(z') + \langle \nabla g(z'), z - z' \rangle + \frac{1}{2}\|z - z'\|^2_{\widehat{\Sigma}_g}. \tag{1.6}$$

For any given $(y', z') \in \mathcal{Y} \times \mathcal{Z}$, $\sigma \in (0, +\infty)$ and any $(y, z, x) \in \mathcal{Y} \times \mathcal{Z} \times \mathcal{X}$, the majorized augmented Lagrangian function is defined as

$$\begin{aligned}\widehat{\mathcal{L}}_\sigma(y, z, x; y', z') := \; & \hat{f}(y, y') + p(y) + \hat{g}(z, z') + q(z) + \langle x, \mathcal{A}^* y + \mathcal{B}^* z - c \rangle \\ & + \frac{\sigma}{2}\|\mathcal{A}^* y + \mathcal{B}^* z - c\|^2,\end{aligned} \tag{1.7}$$

where $\hat{f}$ and $\hat{g}$ are defined by (1.5) and (1.6), respectively. Then the majorized iPADMM can be described as follows.

---
**Algorithm 1** A majorized iPADMM

---
Let $\sigma > 0$ and $\tau \in (0, (1+\sqrt{5})/2)$ be given parameters. Let $\mathcal{S}$ and $\mathcal{T}$ be given self-adjoint, possibly indefinite, linear operators defined on $\mathcal{Y}$ and $\mathcal{Z}$, respectively. Input $(y^0, z^0, x^0) \in \operatorname{dom} p \times \operatorname{dom} q \times \mathcal{X}$. Set $k := 0$.

**Step 1.** Compute

$$\begin{cases} y^{k+1} \in \arg\min_{y \in \mathcal{Y}} \widehat{\mathcal{L}}_\sigma(y, z^k, x^k; y^k, z^k) + \frac{1}{2}\|y - y^k\|^2_\mathcal{S} & (1.8a) \\ \quad = \arg\min_{y \in \mathcal{Y}} \; p(y) + \langle y, \nabla f(y^k) + \mathcal{A} x^k \rangle + \frac{\sigma}{2}\|\mathcal{A}^* y + \mathcal{B}^* z^k - c\|^2 + \frac{1}{2}\|y - y^k\|^2_{\widehat{\Sigma}_f + \mathcal{S}}, \\ z^{k+1} \in \arg\min_{z \in \mathcal{Z}} \widehat{\mathcal{L}}_\sigma(y^{k+1}, z, x^k; y^k, z^k) + \frac{1}{2}\|z - z^k\|^2_\mathcal{T} & (1.8b) \\ \quad = \arg\min_{z \in \mathcal{Z}} \; q(z) + \langle z, \nabla g(z^k) + \mathcal{B} x^k \rangle + \frac{\sigma}{2}\|\mathcal{A}^* y^{k+1} + \mathcal{B}^* z - c\|^2 + \frac{1}{2}\|z - z^k\|^2_{\widehat{\Sigma}_g + \mathcal{T}}, \\ x^{k+1} := x^k + \tau\sigma(\mathcal{A}^* y^{k+1} + \mathcal{B}^* z^{k+1} - c). & (1.8c) \end{cases}$$

**Step 2.** If a termination criterion is not met, set $k := k + 1$ and go to Step 1.

---

As mentioned in Han et al. [27], a large amount of literature focus on designing variant forms of ADMM and their applications. The literature on linear rate convergence, especially the Q-linear rate, on the other hand, is relatively sparse. Based on the close connections among DR, PPA and classic ADMM with $\tau = 1$, the R-linear rate convergence of ADMM can be established under conditions such as strong monotonicity, error bounds and etc. By using some different technical tools, the local Q-linear rate convergence of sequence $\{(z^k, x^k)\}$, generated by the ADMM as defined in Eckstein and Berstekas [15], is established by Han and Yuan [28]. Yang and Han [49] show the Q-linear rate convergence of ADMM for solving a class of convex piecewise linear-quadratic problems. In the same paper, the Q-linear rate convergence of linearized ADMM (L-ADMM) with positive definite proximal operators and $\tau = 1$ is also proved. Under assumptions that either $f(y) + p(y)$ or $g(z) + q(z)$ is strongly convex with a Lipschitz continuous gradient and others, Deng and Yin [9] provide some scenarios on the R-linear and Q-linear rates of convergence of ADMM. Recently, based on the easy-to-use convergence theorem in Fazel et al. [17], the Q-linear rate convergence theorem for sPADMM with $\tau \in (0, (1 + \sqrt{5})/2)$ is established in [27] under a calmness condition, which holds automatically for convex composite piecewise-linear programming. To know more about the convergence rate analysis of ADMM, we refer to [27] and the references therein. Inspired by [27], we resolve the Q-linear rate convergence issue for the majorized iPADMM scheme with $\tau \in (0, (1 + \sqrt{5})/2)$, in addition to some improvements on the global convergence analysis.

The main contributions of this paper are threefold. First, we refine the conditions in [34, Theorem 10 (b)] and establish the global convergence of the majorized iPADMM with

$$\frac{1}{2}\Sigma_f + \mathcal{S} + \sigma\mathcal{A}\mathcal{A}^* \succ 0 \ \& \ \frac{1}{2}\Sigma_g + \mathcal{T} + \sigma\mathcal{B}\mathcal{B}^* \succ 0 \tag{1.9}$$

being replaced by

$$\frac{1}{2}\widehat{\Sigma}_f + \mathcal{S} + \sigma\mathcal{A}\mathcal{A}^* \succ 0 \ \& \ \frac{1}{2}\widehat{\Sigma}_g + \mathcal{T} + \sigma\mathcal{B}\mathcal{B}^* \succ 0. \tag{1.10}$$

Note that conditions (1.9) coincide with (1.10) when $f$ and $g$ are linear/quadratic. However, the improvement from (1.10) is significant when $f$ or $g$ is not linear/quadratic. In fact, there are many widely used loss functions in statistical inference and machine learning that are not linear/quadratic, such as the logistic loss function, the multinomial logistic loss function (see e.g. [19]), and the cox loss function (see e.g. [47]).

Second, we build up a general Q-linear rate convergence theorem based on an inequality associated with the iteration sequences generated by majorized iPADMM. Note that the convergence results on sPADMM [27] are no longer applicable due to the indefiniteness of the proximal terms.

Third, inspired by the global convergence and excellent numerical performance of symmetric Gaussian-Seidel based multi-block sPADMM (sGS-sPADMM) for multi-block linearly constrained convex programming, we present a symmetric Gaussian-Seidel based multi-block majorized iPADMM (majorized sGS-iPADMM). The linear rate convergence results for majorized sGS-iPADMM are established through converting it into an equivalent 2-block majorized iPADMM. This is one of the most important motivations that we consider the objective function in (1.1) in the form of $f(y) + p(y)$ and $g(z) + q(z)$. It is worth mentioning that the sGS-sPDMM was initially presented by Li et al. [35]. For more discussions about symmetric Gaussian-Seidel techniques, we refer to [3, 36, 35, 37] and the references therein.

The rest of the paper is organized as follows. In section 2, we give some preliminaries that will be frequently used in other sections. In section 3, we refine the convergence result [34, Theorem 10 (b)] and then establish a general Q-linear rate convergence theorem under a metric subregularity condition. In section 4, we propose a majorized sGS-iPADMM for the multi-block composite optimization problem. Moreover, we show its convergence rate by establishing the relationship between the majorized sGS-iPADMM and the 2-block majorized iPADMM. In section 5, we apply the majorized iPADMM to three types of regularized logistic regression and then present the numerical results. Finally, we give some concluding remarks and future works in section 6.

## 2 Preliminaries

In this section, we summarize and study some preliminaries that will be used in the subsequent analysis. Let $(\bar{y}, \bar{z})$ be the optimal solution of problem (1.1). If there exists $\bar{x} \in \mathcal{X}$ such that $(\bar{y}, \bar{z}, \bar{x})$ satisfies the following Karush-Kuhn-Tucker (KKT) system:

$$\begin{cases} 0 \in \partial p(y) + \nabla f(y) + \mathcal{A}x, \\ 0 \in \partial q(z) + \nabla g(z) + \mathcal{B}x, \\ \mathcal{A}^*y + \mathcal{B}^*z - c = 0, \end{cases} \quad (2.1)$$

then, $(\bar{y}, \bar{z}, \bar{x})$ is called a KKT point of problem (1.1). Denote the set of KKT points by $\overline{\Omega}$.

For any convex function $\theta : \mathcal{Y} \to \Re \cup \{+\infty\}$, the Moreau-Yosida proximal mapping $\mathrm{Pr}_\theta(\cdot)$ associated with $\theta$ is defined by

$$\mathrm{Pr}_\theta(y) := \arg\min_{y' \in \mathcal{Y}} \left\{ \theta(y') + \frac{1}{2}\|y' - y\|^2 \right\}, \; \forall y \in \mathcal{Y}.$$

It is well known that the Moreau-Yosida proximal mapping $\mathrm{Pr}_\theta(\cdot)$ is globally Lipschitz continuous with modulus one, see e.g. [29, 33].

Denote $u := (y, z, x) \in \mathcal{U}$ with $\mathcal{U} := \mathcal{Y} \times \mathcal{Z} \times \mathcal{X}$. Define the KKT mapping $\mathcal{R} : \mathcal{U} \to \mathcal{U}$ as

$$\mathcal{R}(u) := \begin{pmatrix} y - \mathrm{Pr}_p[y - (\nabla f(y) + \mathcal{A}x)] \\ z - \mathrm{Pr}_q[z - (\nabla g(z) + \mathcal{B}x)] \\ c - \mathcal{A}^*y - \mathcal{B}^*z \end{pmatrix}. \quad (2.2)$$

From [40], we know that $u \in \overline{\Omega}$ if and only if $\mathcal{R}(u) = 0$.

Let $\mathcal{F} : \mathcal{X} \rightrightarrows \mathcal{Y}$ be a multivalued mapping. Denote its inverse by $\mathcal{F}^{-1}$. Define the graph of multivalued function $\mathcal{F}$ as follows

$$\mathrm{gph}\,\mathcal{F} := \{(x, y) \in \mathcal{X} \times \mathcal{Y} \mid y \in \mathcal{F}(x)\}.$$

**Definition 2.1.** *A multivalued mapping $\mathcal{F} : \mathcal{X} \rightrightarrows \mathcal{Y}$ is said to be metrically subregular at $(\bar{x}, \bar{y}) \in \mathrm{gph}\,\mathcal{F}$ with modulus $\eta > 0$ if there exists a neighborhood $\mathcal{U}$ of $\bar{x}$ such that*

$$\mathrm{dist}(x, \mathcal{F}^{-1}(\bar{y})) \leq \eta \mathrm{dist}(\bar{y}, \mathcal{F}(x)), \; \forall x \in \mathcal{U}.$$

The definition of metric subregularity is directly from [13, Definition 3.1]. It is well known that $\mathcal{F}$ is metrically subregular at $(\bar{x}, \bar{y}) \in \mathrm{gph}\,\mathcal{F}$ if and only if its inverse multivalued mapping $\mathcal{F}^{-1}$ is calm (c.f. [50, Definition 2.6], [12, 3.8(3H)]) at $(\bar{y}, \bar{x}) \in \mathrm{gph}\,\mathcal{F}^{-1}$.

From [42, Proposition 1] and [45], we know that if $\mathcal{F}$ is piecewise polyhedral or $\mathcal{F}$ is the subdifferential mapping of a convex piecewise linear-quadratic function, then $\mathcal{F}$ is metrically subregular at $(\bar{x}, \bar{y}) \in \mathrm{gph}\,\mathcal{F}$. Till now, numerous works have been done to study the sufficient conditions of calmness of KKT solution mappings, we refer to [10, 8, 36, 25] and the references therein.

In order to establish the linear rate convergence of the majorized iPADMM, we need the metric subregularity of the KKT mapping $\mathcal{R}$. From the Definition 2.1, the metric subregularity of $\mathcal{R}$ at $(\bar{u}, 0) \in \mathrm{gph}\,\mathcal{R}$ with modulus $\eta > 0$ can be described as: there exists a scalar $\rho > 0$ such that

$$\mathrm{dist}(u, \overline{\Omega}) \leq \eta \|\mathcal{R}(u)\|, \; \forall u \in \{u \in \mathcal{U} : \|u - \bar{u}\| \leq \rho\}. \quad (2.3)$$

The above condition is also referred to as the existence of a local error bound. From [36, Theorem 1 & Remark 1], we know that the KKT mappings corresponding to the Lasso, elastic net Lasso regularized logistic regression models are metrically subregular at $(\bar{u}, 0) \in \mathrm{gph}\,\mathcal{R}$.

Since for any proper closed convex function $\theta : \mathcal{Y} \to \Re \cup \{+\infty\}$, the subdifferential $\partial \theta(\cdot)$ is a monotone multivalued mapping (see [44]), i.e., for any $y^1, y^2 \in \text{dom}\theta$, it holds that

$$\langle \zeta^1 - \zeta^2, y^1 - y^2 \rangle \geq 0, \ \forall \zeta^1 \in \partial \theta(y^1), \forall \zeta^2 \in \partial \theta(y^2). \tag{2.4}$$

**Proposition 2.1.** *[5, Proposition 2.6.5] Let multivalued mapping $\mathcal{F}$ be Lipschitz on an open convex set $\mathcal{V}$ in $\Re^n$, and let $x$ and $y$ be points in $\mathcal{V}$. Then one has*

$$\mathcal{F}(y) - \mathcal{F}(x) \in \text{conv}\, \partial \mathcal{F}([x,y])(y-x),$$

*where $\text{conv}\, \partial \mathcal{F}([x,y])(y-x)$ denotes the convex hull of all points of the form $\zeta(y-x)$, where $\zeta \in \partial \mathcal{F}(u)$ for some point $u$ in $[x,y]$.*

Throughout the subsequent analysis, we always assume that the following two assumptions hold.

**Assumption 2.1.** *The KKT system (2.1) has at least one solution, i.e., $\overline{\Omega} \neq \emptyset$.*

**Assumption 2.2.** *The two self-adjoint linear operators $\mathcal{S} : \mathcal{Y} \to \mathcal{Y}$ and $\mathcal{T} : \mathcal{Z} \to \mathcal{Z}$ in majorized iPADMM satisfy*

$$\mathcal{S} \succeq -\frac{1}{2}\widehat{\Sigma}_f \ \& \ \mathcal{T} \succeq -\frac{1}{2}\widehat{\Sigma}_g. \tag{2.5}$$

**Remark 2.1.** *Assumption 2.2 means that the proximal terms $\mathcal{S}$ and $\mathcal{T}$ cannot be too indefinite as long as $\widehat{\Sigma}_f$ and $\widehat{\Sigma}_g$ are not very big. Note that $\widehat{\Sigma}_f$ and $\widehat{\Sigma}_g$ should be chosen as small as possible provided (1.5) and (1.6) are satisfied. For example, when $f$ is a convex quadratic function, we choose $\widehat{\Sigma}_f = \Sigma_f = \nabla^2 f$, where $\nabla^2 f$ is the Hessian matrix of $f$.*

## 3 Q-Linear rate of convergence of the majorized iPADMM

This section aims to analyze the convergence rate of the majorized iPADMM for solving (1.1). We show that the algorithm achieves a Q-linear rate of convergence under some mild conditions. Before formally stating our main results, we first give some technical results.

### 3.1 Technical lemmas

For notational convenience, for any $\tau \in (0, +\infty)$, define

$$s_\tau := \frac{5 - \tau - 3\min(\tau, \tau^{-1})}{4}, \ t_\tau := \frac{1 - \tau + \min(\tau, \tau^{-1})}{2},$$

and two self-adjoint linear operators:

$$\mathcal{M} := \text{Diag}\left(\widehat{\Sigma}_f + \mathcal{S}, \widehat{\Sigma}_g + \mathcal{T} + \sigma \mathcal{B} \mathcal{B}^*, (\tau\sigma)^{-1}\mathcal{I}\right) + s_\tau \sigma \mathcal{E} \mathcal{E}^*, \tag{3.1}$$

$$\mathcal{H} := \text{Diag}\left(\mathcal{H}_f, \mathcal{H}_g, t_\tau(\tau^2\sigma)^{-1}\mathcal{I}\right) + \frac{1}{4}t_\tau \sigma \mathcal{E} \mathcal{E}^*, \tag{3.2}$$

where

$$\mathcal{H}_f := \frac{1}{2}\widehat{\Sigma}_f + \mathcal{S}, \ \mathcal{H}_g := \frac{1}{2}\widehat{\Sigma}_g + \mathcal{T} + 2t_\tau \tau \sigma \mathcal{B} \mathcal{B}^*, \tag{3.3}$$

and $\mathcal{E} : \mathcal{X} \to \mathcal{Y} \times \mathcal{Z} \times \mathcal{X}$ is the linear operator such that its adjoint $\mathcal{E}^*$ satisfies $\mathcal{E}^*(y, z, x) = \mathcal{A}^* y + \mathcal{B}^* z$. For a given self-adjoint linear operator $\mathcal{G} : \mathcal{X} \to \mathcal{X}$, we denote the largest eigenvalue by $\lambda_{\max}(\mathcal{G})$ and for any $k \geq 0$,
$$r^k := \mathcal{A}^* y^k + \mathcal{B}^* z^k - c.$$

The proofs of the lemmas in this subsection are all presented in Appendix for readability.

**Lemma 3.1.** *Let $\{u^k := (y^k, z^k, x^k)\}$ be the infinite sequence generated by the majorized iPADMM. Then for any $k \geq 0$,*
$$\|\mathcal{R}(u^{k+1})\|^2 \leq \|u^{k+1} - u^k\|^2_{\mathcal{H}_0}, \tag{3.4}$$
*where*
$$\mathcal{H}_0 := \max\{\kappa_1, \kappa_2, \kappa_3\} \mathrm{Diag}\left(\mathcal{I}, \mathcal{I} + \sigma \mathcal{B}\mathcal{B}^*, (\tau^2 \sigma)^{-1} \mathcal{I}\right)$$
*with*
$$\kappa_1 := 3(\lambda_{\max}(\mathcal{S} + \tfrac{1}{2}\widehat{\Sigma}_f) + \tfrac{1}{2}\lambda_{\max}(\widehat{\Sigma}_f))^2,$$
$$\kappa_2 := \max\{2(\lambda_{\max}(\mathcal{T} + \tfrac{1}{2}\widehat{\Sigma}_g) + \tfrac{1}{2}\lambda_{\max}(\widehat{\Sigma}_g))^2, 3\sigma \lambda_{\max}(\mathcal{A}^* \mathcal{A})\},$$
$$\kappa_3 := \sigma^{-1} + (1-\tau)^2 \sigma(3\lambda_{\max}(\mathcal{A}^* \mathcal{A}) + 2\lambda_{\max}(\mathcal{B}^* \mathcal{B})).$$

The above lemma is inspired by [27, Lemma 1], but its proof is more complicated due to the majorization techniques. In order to refine the global convergence results of the majorized iPADMM presented in [34, Theorem 10(b)], we need the following lemma.

**Lemma 3.2.** *Let $h : \mathcal{X} \to \Re$ be a smooth convex function and there is a self-adjoint positive semidefinite linear operator $\mathcal{P}$ such that, for any given $\bar{x} \in \mathcal{X}$,*
$$h(x) \leq h(\bar{x}) + \langle \nabla h(\bar{x}), x - \bar{x} \rangle + \frac{1}{2}\|x - \bar{x}\|^2_{\mathcal{P}}, \quad \forall x \in \mathcal{X}. \tag{3.5}$$
*Then it holds that*
$$\langle \nabla h(x) - \nabla h(\bar{x}), y - \bar{x} \rangle \geq -\frac{1}{4}\|x - y\|^2_{\mathcal{P}}, \quad \forall x, y \in \mathcal{X}. \tag{3.6}$$

**Lemma 3.3.** *Let $\{(y^k, z^k, x^k)\}$ be the infinite sequence generated by the majorized iPADMM. Then, for any $\bar{u} := (\bar{y}, \bar{z}, \bar{x}) \in \overline{\Omega}$, $\tau > 0$ and $k \geq 0$, we have*
$$\begin{aligned}\phi_k - \phi_{k+1} \geq\ & \|y^{k+1} - y^k\|^2_{\frac{1}{2}\widehat{\Sigma}_f + \mathcal{S}} + \|z^{k+1} - z^k\|^2_{\frac{1}{2}\widehat{\Sigma}_g + \mathcal{T}} \\ & + (1-\tau)\sigma \|r^{k+1}\|^2 + \sigma \|\mathcal{A}^* y^{k+1} + \mathcal{B}^* z^k - c\|^2,\end{aligned} \tag{3.7}$$
*where for any $k \geq 0$,*
$$\phi_k := (\tau\sigma)^{-1}\|x^k - \bar{x}\|^2 + \|y^k - \bar{y}\|^2_{\widehat{\Sigma}_f + \mathcal{S}} + \|z^k - \bar{z}\|^2_{\widehat{\Sigma}_g + \mathcal{T} + \sigma \mathcal{B}\mathcal{B}^*}. \tag{3.8}$$

Since the proof of the following lemma is not much different from the one in [34, Theorem 10, Inequality (55)] except for replacing Inequality (33) in [34] by (3.7) in Lemma 3.3, we include an outline in Appendix A.4.

**Lemma 3.4.** *Let $\{u^k := (y^k, z^k, x^k)\}$ be the infinite sequence generated by the majorized iPADMM. For each $k$ and any KKT point $\bar{u} := (\bar{y}, \bar{z}, \bar{x})$, let $\phi_k$ be defined in (3.8). Then, for any $k \geq 1$, one has*

$$\begin{aligned}&\left[\phi_k + (1 - \min(\tau, \tau^{-1}))\sigma \|r^k\|^2 + \|z^k - z^{k-1}\|^2_{\widehat{\Sigma}_g + \mathcal{T}}\right] \\ &- \left[\phi_{k+1} + (1 - \min(\tau, \tau^{-1}))\sigma \|r^{k+1}\|^2 + \|z^{k+1} - z^k\|^2_{\widehat{\Sigma}_g + \mathcal{T}}\right] \\ &\geq t_{k+1} + (-\tau + \min(1 + \tau, 1 + \tau^{-1}))\sigma \|r^{k+1}\|^2,\end{aligned} \quad (3.9)$$

*where*

$$t_{k+1} := \|y^{k+1} - y^k\|^2_{\mathcal{H}_f} + \|z^{k+1} - z^k\|^2_{\mathcal{H}_g}. \quad (3.10)$$

The following lemma plays an important role in establishing the global convergence theorem for the majorized iPADMM

**Lemma 3.5.** *Let $\tau \in (0, (1+\sqrt{5})/2)$, $\mathcal{M}$ and $\mathcal{H}$ be defined in (3.1) and (3.2), respectively. Then,*

$$\frac{1}{2}\widehat{\Sigma}_f + \mathcal{S} + \sigma \mathcal{A}\mathcal{A}^* \succ 0 \ \& \ \frac{1}{2}\widehat{\Sigma}_g + \mathcal{T} + \sigma \mathcal{B}\mathcal{B}^* \succ 0 \Leftrightarrow \mathcal{H} \succ 0 \Rightarrow \mathcal{M} \succ 0. \quad (3.11)$$

## 3.2 Convergence analysis

In this subsection, we investigate the rate of convergence of majorized iPADMM for solving (1.1). Inspired by [27, Proposition 4], we first develop a key inequality needed for proving the linear rate convergence for the majorized iPADMM.

**Proposition 3.1.** *Let $\tau \in (0, (1+\sqrt{5})/2)$ and $\{u^k := (y^k, z^k, x^k)\}$ be the infinite sequence generated by the majorized iPADMM. Then for any KKT point $\bar{u} := (\bar{y}, \bar{z}, \bar{x})$ and any $k \geq 1$,*

$$\begin{aligned}&\|u^{k+1} - \bar{u}\|^2_{\mathcal{M}} + \|z^{k+1} - z^k\|^2_{\widehat{\Sigma}_g + \mathcal{T}} \\ &\leq \|u^k - \bar{u}\|^2_{\mathcal{M}} + \|z^k - z^{k-1}\|^2_{\widehat{\Sigma}_g + \mathcal{T}} - \|u^{k+1} - u^k\|^2_{\mathcal{H}}.\end{aligned} \quad (3.12)$$

*Consequently, we have for all $k \geq 1$,*

$$\begin{aligned}&\operatorname{dist}^2_{\mathcal{M}}(u^{k+1}, \overline{\Omega}) + \|z^{k+1} - z^k\|^2_{\widehat{\Sigma}_g + \mathcal{T}} \\ &\leq \operatorname{dist}^2_{\mathcal{M}}(u^k, \overline{\Omega}) + \|z^k - z^{k-1}\|^2_{\widehat{\Sigma}_g + \mathcal{T}} - \|u^{k+1} - u^k\|^2_{\mathcal{H}}.\end{aligned} \quad (3.13)$$

*Proof.* By reorganizing the inequality in (3.9), one has

$$\begin{aligned}&\left[\phi_{k+1} + s_\tau \sigma \|r^{k+1}\|^2 + \|z^{k+1} - z^k\|^2_{\widehat{\Sigma}_g + \mathcal{T}}\right] - \left[\phi_k + s_\tau \sigma \|r^k\|^2 + \|z^k - z^{k-1}\|^2_{\widehat{\Sigma}_g + \mathcal{T}}\right] \\ &\leq -\left\{t_\tau \sigma \|r^{k+1}\|^2 + \tfrac{1}{2}\sigma t_\tau (\|r^{k+1}\|^2 + \|r^k\|^2) + t_{k+1}\right\}.\end{aligned} \quad (3.14)$$

From definitions of $x^{k+1}$ and $(\bar{y}, \bar{z}, \bar{x})$, we have

$$\begin{aligned}r^{k+1} &= (\tau\sigma)^{-1}(x^{k+1} - x^k) = \mathcal{A}^*(y^{k+1} - \bar{y}) + \mathcal{B}^*(z^{k+1} - \bar{z}), \\ r^k &= \mathcal{A}^*(y^k - \bar{y}) + \mathcal{B}^*(z^k - \bar{z}), \\ \|r^{k+1}\|^2 + \|r^k\|^2 &\geq \tfrac{1}{2}\|\mathcal{A}^*(y^{k+1} - y^k) + \mathcal{B}^*(z^{k+1} - z^k)\|^2.\end{aligned} \quad (3.15)$$

Then we can get (3.12) by substituting (3.15) into (3.14). Since (3.12) holds for any $\bar{u} \in \overline{\Omega}$, we can get (3.13) from the fact that $\overline{\Omega}$ is a nonempty closed convex set. The proof is completed. □

Now we are ready to establish the global convergence and the linear rate of convergence for the majorized iPADMM under a metric subregularity condition of $\mathcal{R}$ at some $(\bar{u}, 0) \in \text{gph}\mathcal{R}$.

**Theorem 3.1.** *Let $\tau \in (0, (1+\sqrt{5})/2)$, the two self-adjoint linear operators $\mathcal{S}$ and $\mathcal{T}$ satisfy (2.5),*

$$\frac{1}{2}\widehat{\Sigma}_f + \mathcal{S} + \sigma\mathcal{A}\mathcal{A}^* \succ 0 \,\&\, \frac{1}{2}\widehat{\Sigma}_g + \mathcal{T} + \sigma\mathcal{B}\mathcal{B}^* \succ 0, \tag{3.16}$$

*and $\{u^k := (y^k, z^k, x^k)\}$ be the infinite sequence generated by the majorized iPADMM. Then, one has the following results.*

(a) *The sequence $\{(y^k, z^k)\}$ converges to an optimal solution of problem (1.1) and $\{x^k\}$ converges to an optimal solution of the dual of problem (1.1).*

(b) *Suppose that the sequence $\{(y^k, z^k, x^k)\}$ converges to a KKT point $\bar{u} := (\bar{y}, \bar{z}, \bar{x})$ and the KKT mapping $\mathcal{R}$ is metrically subregular at $(\bar{u}, 0) \in \text{gph}\mathcal{R}$ with modulus $\eta > 0$. Then there exist a positive number $\mu \in (0, 1)$ and an integer $k_0 \geq 1$ such that for all $k \geq k_0$,*

$$\text{dist}_{\mathcal{M}}^2(u^{k+1}, \overline{\Omega}) + \|z^{k+1} - z^k\|_{\widehat{\Sigma}_g + \mathcal{T}}^2 \leq \mu \left[ \text{dist}_{\mathcal{M}}^2(u^k, \overline{\Omega}) + \|z^k - z^{k-1}\|_{\widehat{\Sigma}_g + \mathcal{T}}^2 \right], \tag{3.17}$$

*Moreover, there exists a positive number $\hat{\mu} \in [\mu, 1)$ such that for all $k \geq 1$,*

$$\text{dist}_{\mathcal{M}}^2(u^{k+1}, \overline{\Omega}) + \|z^{k+1} - z^k\|_{\widehat{\Sigma}_g + \mathcal{T}}^2 \leq \hat{\mu} \left[ \text{dist}_{\mathcal{M}}^2(u^k, \overline{\Omega}) + \|z^k - z^{k-1}\|_{\widehat{\Sigma}_g + \mathcal{T}}^2 \right]. \tag{3.18}$$

*Proof.* We first prove the convergence on the sequences $\{(y^k, z^k)\}$ and $\{x^k\}$. From Lemma 3.5 and (3.16), we know that $\mathcal{H} \succ 0$, $\mathcal{M} \succ 0$ and $\widehat{\Sigma}_g + \mathcal{T} \succeq 0$. Then it holds from Proposition 3.1 that $\{u^{k+1}\}$ is bounded and

$$\lim_{k \to \infty} \|u^{k+1} - u^k\| = 0.$$

Consequently, there is a subsequence $\{u^{k_i}\}$ which converges to a cluster point $u^\infty$. From Lemma 3.1, we know that

$$\|\mathcal{R}(u^{k_i})\|^2 \leq \|u^{k_i} - u^{k_i - 1}\|_{\mathcal{H}_0}^2,$$

where $\mathcal{H}_0 \succ 0$. Taking limits on both sides of the above inequality, we obtain $\|\mathcal{R}(u^\infty)\| = 0$. Thus, the subsequence $\{u^{k_i}\}$ converges to $u^\infty \in \overline{\Omega}$. Therefore, the sequence $\{\|u^{k_i+1} - u^\infty\|_{\mathcal{M}}^2 + \|z^{k_i+1} - z^{k_i}\|_{\widehat{\Sigma}_g + \mathcal{T}}^2\}$ converges to 0 as $k_i \to \infty$. Since the subsequence is non-increasing and $\|u^{k+1} - u^k\| \to 0$, we have

$$\lim_{k \to \infty} \|u^k - u^\infty\| = 0.$$

Therefore, the whole sequence $\{u^k\}$ converges to $u^\infty$. This completes the proof of the result $(a)$. Next, we prove $(b)$. From $(a)$, we know that the sequence $\{(y^k, z^k, x^k)\}$ generated by the majorized iPADMM converges to a KKT point $\bar{u} = (\bar{y}, \bar{z}, \bar{x})$. Then there exist $\rho > 0$ and an integer $k_0 \geq 1$ such that for all $k \geq k_0$,

$$\|u^{k+1} - \bar{u}\| \leq \rho.$$

Therefore, by using Lemma 3.1 and (2.3), we know that for all $k \geq k_0$,

$$\text{dist}^2(u^{k+1}, \overline{\Omega}) \leq \eta^2 \|\mathcal{R}(u^{k+1})\|^2 \leq \eta^2 \|u^{k+1} - u^k\|_{\mathcal{H}_0}^2. \tag{3.19}$$

The definition of $\mathcal{H}$ and the fact that $\mathcal{H} \succ 0$ imply that $\mathcal{H}_g \succ 0$. Then there exists a finite real number $\varrho_1 > 0$ such that $\widehat{\Sigma}_g + \mathcal{T} \preceq \varrho_1 \mathcal{H}_g$ and consequently, it holds that

$$\|z^{k+1} - z^k\|^2_{\widehat{\Sigma}_g + \mathcal{T}} \leq \varrho_1 \|u^{k+1} - u^k\|^2_{\mathcal{H}}.$$

Similarly, there exists a finite real number $\varrho_2 > 0$ such that $\mathcal{H}_0 \preceq \varrho_2 \mathcal{H}$. It follows from (3.19) that for all $k \geq k_0$,

$$\begin{aligned} \|u^{k+1} - u^k\|^2_{\mathcal{H}} &\geq \varrho_2^{-1} \|u^{k+1} - u^k\|^2_{\mathcal{H}_0} \\ &\geq \varrho_2^{-1} \eta^{-2} \mathrm{dist}^2_{\mathcal{M}}(u^{k+1}, \overline{\Omega}) \geq \varrho_2^{-1} \eta^{-2} \lambda^{-1}_{\max}(\mathcal{M}) \mathrm{dist}^2_{\mathcal{M}}(u^{k+1}, \overline{\Omega}). \end{aligned} \tag{3.20}$$

Let $\kappa := \frac{1}{1+\varrho_1 \beta}$ with $\beta := \varrho_2^{-1} \eta^{-2} \lambda^{-1}_{\max}(\mathcal{M})$. From (3.13), we have that, for all $k \geq k_0$,

$$\begin{aligned} &\left[\mathrm{dist}^2_{\mathcal{M}}(u^{k+1}, \overline{\Omega}) + \|z^{k+1} - z^k\|^2_{\widehat{\Sigma}_g + \mathcal{T}}\right] - \left[\mathrm{dist}^2_{\mathcal{M}}(u^k, \overline{\Omega}) + \|z^k - z^{k-1}\|^2_{\widehat{\Sigma}_g + \mathcal{T}}\right] \\ &\leq -\left((1-\kappa)\|u^{k+1} - u^k\|^2_{\mathcal{H}} + \kappa \|u^{k+1} - u^k\|^2_{\mathcal{H}}\right) \\ &\leq -\left((1-\kappa)\varrho_1^{-1} \|z^{k+1} - z^k\|^2_{\widehat{\Sigma}_g + \mathcal{T}} + \kappa \varrho_2^{-1} \eta^{-2} \lambda^{-1}_{\max}(\mathcal{M}) \mathrm{dist}^2_{\mathcal{M}}(u^{k+1}, \overline{\Omega})\right). \end{aligned} \tag{3.21}$$

Then by reorganizing the above inequality, we know that for all $k \geq k_0$,

$$\left(1 + \kappa \varrho_2^{-1} \eta^{-2} \lambda^{-1}_{\max}(\mathcal{M})\right) \mathrm{dist}^2_{\mathcal{M}}(u^{k+1}, \overline{\Omega}) + \left(1 + (1-\kappa)\varrho_1^{-1}\right) \|z^{k+1} - z^k\|^2_{\widehat{\Sigma}_g + \mathcal{T}}$$

$$\leq \mathrm{dist}^2_{\mathcal{M}}(u^k, \overline{\Omega}) + \|z^k - z^{k-1}\|^2_{\widehat{\Sigma}_g + \mathcal{T}}.$$

It is easy to check that

$$1 + \kappa \varrho_2^{-1} \eta^{-2} \lambda^{-1}_{\max}(\mathcal{M}) = 1/\mu,$$

where

$$\mu := \frac{\varrho_1 \beta + 1}{1 + \beta + \varrho_1 \beta} < 1. \tag{3.22}$$

Then we know that inequality (3.17) holds.

By combining (3.17) with Lemma 3.1, (3.13) in Proposition 3.1, we can obtain directly that there exists a positive number $\hat{\mu} \in [\mu, 1)$ such that (3.18) holds for all $k \geq 1$. This completes the proof. $\square$

**Remark 3.1.** *We make the following comments.*

(a) *From the result (a) in Theorem 3.1, we can see that the majorized iPADMM is still globally convergent when the conditions $\frac{1}{2}\Sigma_f + \mathcal{S} + \sigma \mathcal{A}\mathcal{A}^* \succ 0$ & $\frac{1}{2}\Sigma_g + \mathcal{T} + \sigma \mathcal{B}\mathcal{B}^* \succ 0$ in [34, Theorem 10 (b)] are replaced by*

$$\frac{1}{2}\widehat{\Sigma}_f + \mathcal{S} + \sigma \mathcal{A}\mathcal{A}^* \succ 0 \ \& \ \frac{1}{2}\widehat{\Sigma}_g + \mathcal{T} + \sigma \mathcal{B}\mathcal{B}^* \succ 0.$$

*As mentioned in the introduction, the improvement is significant when $f$ or $g$ is not linear/quadratic.*

(b) *By observing the expression of parameter $\mu$ in the proof of result (b) in Theorem 3.1, in order to increase the convergence speed, under the premise of satisfying (2.5) and (3.16), the linear operators $\widehat{\Sigma}_f$ and $\widehat{\Sigma}_g$ should be chosen such that the majorized functions $\hat{f}$ and $\hat{g}$ are as close to $f$ and $g$ as possible, and the proximal terms $\mathcal{S}$ and $\mathcal{T}$ should be chosen as close to $-\frac{1}{2}\widehat{\Sigma}_f$ and $-\frac{1}{2}\widehat{\Sigma}_g$ as possible. For simplicity, we assume that $\mathcal{A}\mathcal{A}^* \succ 0$, then one should choose $\mathcal{S} = -\frac{1}{2}\widehat{\Sigma}_f$. This illustrates the claim mentioned in the introduction that the linear rate convergence results for sPADMM established in [27] are no longer applicable for our majorized iPADMM.*

## 4 Application I: A majorized sGS-iPADMM

Consider the following general multi-block convex composite programming model,

$$\begin{aligned} \min_{y,z} \quad & p(y_1) + f(y_1, \ldots, y_s) + q(z_1) + g(z_1, \ldots, z_t) \\ \text{s.t.} \quad & \mathcal{A}^* y + \mathcal{B}^* z = c, \\ & y \in \mathcal{Y},\ z \in \mathcal{Z}, \end{aligned} \quad (4.1)$$

where $s$ and $t$ are given nonnegative integers, $\mathcal{Y} := \mathcal{Y}_1 \times \ldots \times \mathcal{Y}_s$, $\mathcal{Z} := \mathcal{Z}_1 \times \ldots \times \mathcal{Z}_t$, $f(y_1, \ldots, y_s) := \sum_{i=1}^s f_i(y_i)$, and $g(z_1, \ldots, z_t) := \sum_{j=1}^t g_j(z_j)$.

For $i \in \{1, \ldots, s\}$ and $j \in \{1, \ldots, t\}$, we assume that $f_i : \mathcal{Y}_i \to \Re$ and $g_j : \mathcal{Y}_j \to \Re$ are convex functions with Lipschitz continuous gradients. Then, there exist positive semidefinite operators $\widehat{\Sigma}_{f_i}$ and $\widehat{\Sigma}_{g_j}$ such that for given $y_i' \in \mathcal{Y}_i$, $z_j' \in \mathcal{Z}_j$,

$$\begin{aligned} f_i(y_i) &\le \hat{f}_i(y_i, y_i') := f_i(y_i') + \langle \nabla f_i(y_i'), y_i - y_i'\rangle + \tfrac{1}{2}\|y_i - y_i'\|^2_{\widehat{\Sigma}_{f_i}}, \\ g_j(z_j) &\le \hat{g}_j(z_j, z_j') := g_j(z_j') + \langle \nabla g_j(z_j'), z_j - z_j'\rangle + \tfrac{1}{2}\|z_j - z_j'\|^2_{\widehat{\Sigma}_{g_j}}. \end{aligned}$$

Set

$$\widehat{\Sigma}_f := \mathrm{Diag}(\widehat{\Sigma}_{f_1}, \ldots, \widehat{\Sigma}_{f_s})\ \&\ \widehat{\Sigma}_g := \mathrm{Diag}(\widehat{\Sigma}_{g_1}, \ldots, \widehat{\Sigma}_{g_t}). \quad (4.2)$$

Then it holds that

$$\begin{aligned} f(y) &\le \hat{f}(y, y') := f(y') + \langle \nabla f(y'), y - y'\rangle + \tfrac{1}{2}\|y - y'\|^2_{\widehat{\Sigma}_f}, \\ g(z) &\le \hat{g}(z, z') := g(z') + \langle \nabla g(z'), z - z'\rangle + \tfrac{1}{2}\|z - z'\|^2_{\widehat{\Sigma}_g}. \end{aligned}$$

For any given parameter $\sigma > 0$, the majorized augmented Lagrangian function $\widehat{\mathcal{L}}_\sigma(y, z, x; y', z')$ is defined as (1.7) and the sGS based multi-block majorized ADMM with indefinite proximal terms (majorized sGS-iPADMM) is presented in Algorithm 2.

**Algorithm 2** A majorized sGS-iPADMM

Let $\sigma > 0$ and $\tau \in (0, (1+\sqrt{5})/2)$ be given parameters. For $i \in \{1, \ldots, s\}$ and $j \in \{1, \ldots, t\}$, let $\mathcal{S}_i$, $\mathcal{T}_j$ be given self-adjoint, possibly indefinite, linear operators. Input $(y^0, z^0, x^0) \in \text{dom}\, p \times \text{dom}\, q \times \mathcal{X}$. For $k = 0, 1, 2, \ldots$, perform the following steps.

**Step 1a.** (Backward GS sweep) Compute for $i = s, \ldots, 2$,
$$\bar{y}_i^k = \arg\min_{y_i} \widehat{\mathcal{L}}_\sigma(y_{\leq i-1}^k, y_i, \bar{y}_{\geq i+1}^k, z^k, x^k; y^k, z^k) + \frac{1}{2}\|y_i - y_i^k\|_{\mathcal{S}_i}^2,$$
and $y_1^{k+1} = \arg\min_{y_1} \widehat{\mathcal{L}}_\sigma(y_1, \bar{y}_{\geq 2}^k, z^k, x^k; y^k, z^k) + \frac{1}{2}\|y_1 - y_1^k\|_{\mathcal{S}_1}^2$.

**Step 1b.** (Forward GS sweep) Compute for $i = 2, \ldots, s$,
$$y_i^{k+1} = \arg\min_{y_i} \widehat{\mathcal{L}}_\sigma(y_{\leq i-1}^{k+1}, y_i, \bar{y}_{\geq i+1}^k, z^k, x^k; y^k, z^k) + \frac{1}{2}\|y_i - y_i^k\|_{\mathcal{S}_i}^2.$$

**Step 1c.** (Backward GS sweep) Compute for $j = t, \ldots, 2$,
$$\bar{z}_j^k = \arg\min_{z_j} \widehat{\mathcal{L}}_\sigma(y^{k+1}, z_{\leq j-1}^k, z_i, \bar{z}_{\geq j+1}^k, x^k; y^k, z^k) + \frac{1}{2}\|z_j - z_j^k\|_{\mathcal{T}_j}^2,$$
and $z_1^{k+1} = \arg\min_{z_1} \widehat{\mathcal{L}}_\sigma(y^{k+1}, z_1, \bar{z}_{\geq 2}^k, x^k; y^k, z^k) + \frac{1}{2}\|z_1 - z_1^k\|_{\mathcal{T}_1}^2$.

**Step 1d.** (Forward GS sweep) Compute for $j = 2, \ldots, t$,
$$z_i^{k+1} = \arg\min_{z_j} \widehat{\mathcal{L}}_\sigma(y^{k+1}, z_{\leq j-1}^{k+1}, z_j, \bar{z}_{\geq j+1}^k, x^k; y^k, z^k) + \frac{1}{2}\|z_j - z_j^k\|_{\mathcal{T}_j}^2.$$

**Step 2.** Compute $x^{k+1} = x^k + \tau\sigma(\mathcal{A}^* y^{k+1} + \mathcal{B}^* z^{k+1} - c)$.

---

Denote the following two linear operators,
$$\widetilde{\mathcal{M}} := \frac{1}{2}\widehat{\Sigma}_f + \sigma\mathcal{A}\mathcal{A}^* + \text{Diag}(\mathcal{S}_1, \ldots, \mathcal{S}_s), \quad \widetilde{\mathcal{N}} := \frac{1}{2}\widehat{\Sigma}_g + \sigma\mathcal{B}\mathcal{B}^* + \text{Diag}(\mathcal{T}_1, \ldots, \mathcal{T}_t),$$

where for $i \in \{1, \ldots, s\}$ and $j \in \{1, \ldots, t\}$, $\mathcal{S}_i \succeq -\frac{1}{2}\widehat{\Sigma}_{f_i}$ and $\mathcal{T}_j \succeq -\frac{1}{2}\widehat{\Sigma}_{g_j}$ are self-adjoint linear operators such that the $i$-th diagonal block operator of $\widetilde{\mathcal{M}}$ and $j$-th diagonal block operator of $\widetilde{\mathcal{N}}$ are positive definite, i.e.,
$$\widetilde{\mathcal{M}}_{ii} = \frac{1}{2}\widehat{\Sigma}_{f_i} + \mathcal{S}_i + \sigma\mathcal{A}_i\mathcal{A}_i^* \succ 0 \ \& \ \widetilde{\mathcal{N}}_{jj} = \frac{1}{2}\widehat{\Sigma}_{g_j} + \mathcal{T}_j + \sigma\mathcal{B}_j\mathcal{B}_j^* \succ 0.$$

Moreover, define
$$\mathcal{S} := \text{Diag}(\mathcal{S}_1, \ldots, \mathcal{S}_s) + \text{sGS}(\widetilde{\mathcal{M}}) \ \& \ \mathcal{T} := \text{Diag}(\mathcal{T}_1, \ldots, \mathcal{T}_t) + \text{sGS}(\widetilde{\mathcal{N}}), \tag{4.3}$$

where $\text{sGS}(\widetilde{\mathcal{M}}) := \widetilde{\mathcal{M}}_u \widetilde{\mathcal{M}}_d^{-1} \widetilde{\mathcal{M}}_u^*$ and $\text{sGS}(\widetilde{\mathcal{N}}) := \widetilde{\mathcal{N}}_u \widetilde{\mathcal{N}}_d^{-1} \widetilde{\mathcal{N}}_u^*$. Notations $\widetilde{\mathcal{M}}_u$ ($\widetilde{\mathcal{M}}_d$) and $\widetilde{\mathcal{N}}_u$ ($\widetilde{\mathcal{N}}_d$) stand for the strictly upper triangular (diagonal) block operators of $\widetilde{\mathcal{M}}$ and $\widetilde{\mathcal{N}}$, respectively.

By using the same process in Li et al. [35] (see also Chen et al. [3]), the majorized sGS-iPADMM can be equivalently converted into the following 2-block majorized iPADMM,

$$\begin{cases} y^{k+1} = \arg\min_{y \in \mathcal{Y}} \ p(y) + \langle y, \nabla f(y^k) + \mathcal{A}x^k \rangle + \frac{\sigma}{2}\|\mathcal{A}^*y + \mathcal{B}^*z^k - c\|^2 + \frac{1}{2}\|y - y^k\|^2_{\widehat{\Sigma}_f + \mathcal{S}}, \\ z^{k+1} = \arg\min_{z \in \mathcal{Z}} \ q(z) + \langle z, \nabla g(z^k) + \mathcal{B}x^k \rangle + \frac{\sigma}{2}\|\mathcal{A}^*y^{k+1} + \mathcal{B}^*z - c\|^2 + \frac{1}{2}\|z - z^k\|^2_{\widehat{\Sigma}_g + \mathcal{T}}, \\ x^{k+1} := x^k + \tau\sigma(\mathcal{A}^*y^{k+1} + \mathcal{B}^*z^{k+1} - c), \end{cases}$$

where the operators $\widehat{\Sigma}_f$ and $\widehat{\Sigma}_g$ are defined by (4.2), the proximal terms $\mathcal{S}$ and $\mathcal{T}$ are defined by (4.3). It follows from the choice of $\mathcal{S}_i$ and $\mathcal{T}_j$ that Assumption 2.2 holds and

$$\frac{1}{2}\widehat{\Sigma}_f + \mathcal{S} + \sigma\mathcal{A}\mathcal{A}^* \succ 0 \ \& \ \frac{1}{2}\widehat{\Sigma}_g + \mathcal{T} + \sigma\mathcal{B}\mathcal{B}^* \succ 0.$$

Therefore, directly from Theorem 3.1, we can get the following convergence results for our majorized sGS-iPADMM.

**Proposition 4.1.** *Let $\{u^k := (y^k, z^k, x^k)\}$ be the infinite sequence generated by the majorized sGS-iPADMM and proximal terms $\mathcal{S}$ and $\mathcal{T}$ be defined by (4.3). Then, we have the following results.*

(a) *The sequence $\{(y^k, z^k)\}$ converges to an optimal solution of problem (4.1) and $\{x^k\}$ converges to an optimal solution of the dual of problem (4.1).*

(b) *Suppose that the sequence $\{(y^k, z^k, x^k)\}$ converges to a KKT point $\bar{u} := (\bar{y}, \bar{z}, \bar{x})$ and the KKT mapping $\mathcal{R}$ is metrically subregular at $(\bar{u}, 0) \in \mathrm{gph}\mathcal{R}$. Then the sequence $\{u^k\}$ is linearly convergent to $\bar{u}$.*

Next we use an example to illustrate the fact that the proximal terms generated by (4.3) could be indefinite.

**Example 4.1.** *Consider the following sparse group Lasso (see e.g. [20]),*

$$\min_x \ \frac{1}{2}\|Dx - d\|^2 + \lambda_1\|x\|_1 + \lambda_2 \sum_{l=1}^m w_l\|x_{G_l}\|, \tag{4.4}$$

*where $\lambda_1, \lambda_2 > 0$, $l = 1, 2, \ldots, m, w_l > 0, G_l \subseteq \{1, 2, \ldots, n\}$ contains the indices corresponding to the l-th group of features. For each $l \in \{1, 2, \ldots, g\}$, define the linear operator $\mathcal{P}_l$ by $\mathcal{P}_l x = x_{G_l}$ and $\mathcal{P} = [\mathcal{P}_1; \mathcal{P}_2; \ldots; \mathcal{P}_g]$. The dual formulation (equivalent minimization form) of (4.4) takes the following form:*

$$\begin{aligned} \min_{y,z,\eta} & \ \tfrac{1}{2}\|\theta\|^2 + d^T\theta + \mathbb{I}_{\mathbb{B}_\infty}(\eta) + \mathbb{I}_{\mathbb{B}_2}(z) \\ \mathrm{s.t.} & \ D^T\theta + \eta + \mathcal{P}^*z = 0, \end{aligned} \tag{4.5}$$

*where $\mathbb{B}_\infty := \{\eta \,|\, \|\eta\|_\infty \leq \lambda_1\}$, $\mathbb{B}_2 := \{z \,|\, \|z_{[l]}\|_2 \leq \lambda_2 w_l, \ z_{[l]} \in \Re^{|G_l|}, \ l = 1, \ldots, m\}$ and $\mathcal{P}^*$ is the adjoint operator of $\mathcal{P}$. Let $y := [\theta; \eta]$ and define*

$$f(y) := \frac{1}{2}\|\theta\|^2 + d^T\theta, \ \ p(y) := \mathbb{I}_{\mathbb{B}_\infty}(\eta), \ \ q(z) := \mathbb{I}_{\mathbb{B}_2}(z), \ \ \mathcal{A}^* := [D^T \ I], \ \ \mathcal{B}^* := \mathcal{P}^*.$$

*Consequently, the dual problem (4.5) can be reformulated into the framework of (1.1) and the majorized sGS-iPADMM iterative scheme can be described as below.*

$$\begin{cases} \theta^{k+\frac{1}{2}} = \arg\min_\theta \ \langle \theta, \theta^k + d + Dx^k \rangle + \frac{\sigma}{2}\|D^T\theta + \eta^k + \mathcal{P}^*z^k\|^2 + \frac{1}{2}\|\theta - \theta^k\|^2_{I+\mathcal{S}_1}, \\ \eta^{k+1} = \arg\min_\eta \ \mathbb{I}_{\mathbb{B}_\infty}(\eta) + \frac{\sigma}{2}\|\eta + D^T\theta^{k+\frac{1}{2}} + \mathcal{P}^*z^k + x^k/\sigma\|^2, \\ \theta^{k+1} = \arg\min_\theta \ \langle \theta, \theta^k + d + Dx^k \rangle + \frac{\sigma}{2}\|D^T\theta + \eta^{k+1} + \mathcal{P}^*z^k\|^2 + \frac{1}{2}\|\theta - \theta^k\|^2_{I+\mathcal{S}_1}, \\ z^{k+1} = \arg\min_z \ \mathbb{I}_{\mathbb{B}_2}(z) + \frac{\sigma}{2}\|\mathcal{P}^*z + \eta^{k+1} + D^T\theta^{k+1} + \mathcal{P}x^k/\sigma\|^2, \\ x^{k+1} = x^k + \tau\sigma(D^T\theta^{k+1} + \mathcal{P}^*z^{k+1} + \eta^{k+1}). \end{cases}$$

*Directly from (4.3), we can take*

$$\mathcal{S} = \begin{pmatrix} \mathcal{S}_1 & 0 \\ 0 & 0 \end{pmatrix} + \begin{pmatrix} 0 & 0 \\ \sigma D^T & 0 \end{pmatrix} \begin{pmatrix} \sigma DD^T + \frac{1}{2}I + \mathcal{S}_1 & 0 \\ 0 & \sigma I \end{pmatrix}^{-1} \begin{pmatrix} 0 & \sigma D \\ 0 & 0 \end{pmatrix}, \quad \mathcal{T} = 0.$$

*For simplicity, we assume that $DD^T \succ 0$. Consequently, one can take $\mathcal{S}_1 = -\frac{1}{2}I$ and derive the following indefinite proximal term,*

$$\mathcal{S} = \begin{pmatrix} -\frac{1}{2}I & 0 \\ 0 & \sigma I \end{pmatrix}.$$

## 5 Application II: The regularized logistic regression

In this section, we apply the majorized iPADMM to general regularized logistic regression in the following form,

$$\min_{y,y_0} \ f(y, y_0) + \varphi(y), \tag{5.1}$$

where $f : \Re^{n+1} \to \Re$ is the logistic loss function and $\varphi : \Re^n \to \Re \cup \{+\infty\}$ is a general convex Lasso regularizer. Specifically, the logistic loss function $f$ takes the following form,

$$f(y, y_0) = \frac{1}{N}\sum_{i=1}^N \log(1 + \exp(-b_i(B_i^T y + y_0))),$$

where $B_i \in \Re^n$ are the predictor variables and $b_i \in \{1, -1\}$ are the responses, $i = 1, \ldots, N$. For notation convenience, set $\tilde{y} := [y; y_0] \in \Re^{n+1}$, $A_i := [-b_i B_i; -b_i] \in \Re^{n+1}$ and denote the gradient of $f$ at $\tilde{y} \in \text{dom} f$ by $\nabla f(\tilde{y})$.

Since the gradient $\nabla f$ is Lipschitz continuous on $\text{dom} f$, we know that there exists a positive semidefinite matrix $\widehat{\Sigma}_f$ such that for any given $\tilde{y}' \in \Re^{n+1}$,

$$f(\tilde{y}) \leq \hat{f}(\tilde{y}; \tilde{y}') := f(\tilde{y}') + \langle \nabla f(\tilde{y}'), \tilde{y} - \tilde{y}' \rangle + \frac{1}{2}\|\tilde{y} - \tilde{y}'\|^2_{\widehat{\Sigma}_f}. \tag{5.2}$$

Elementary calculations show that

$$\nabla f(\tilde{y}) = \frac{1}{N}\sum_{i=1}^N \frac{A_i \exp(A_i^T \tilde{y})}{1 + \exp(A_i^T \tilde{y})}, \quad \nabla^2 f(\tilde{y}) = \frac{1}{N}\sum_{i=1}^N A_i A_i^T \frac{\exp(A_i^T \tilde{y})}{(1 + \exp(A_i^T \tilde{y}))^2} \preceq \frac{1}{4N}\sum_{i=1}^N A_i A_i^T.$$

Therefore, the proximal term can be chosen as

$$\widehat{\Sigma}_f := \frac{1}{4N}AA^T, \ A := [A_1, \ldots, A_N] \in \Re^{(n+1)\times N}. \tag{5.3}$$

In this section, we consider the logistic regression with three types of $\varphi(\cdot)$, namely, Lasso, fused Lasso and constrained Lasso. Next, we reformulate these three types of regularized logistic regression into the framework of (1.1) and tailor the majorized iPADMM for each of them. Since our main purpose is to test the numerical performance of our majorized iPADMM, we omit the history and development of these regularized logistic regression models. To know more about these models, we refer to [4, 30, 31, 32, 46, 48, 39].

## 5.1 Fused Lasso logistic regression

The fused Lasso method is introduced by Tibshirani et al. [48] to study the situation that the features have a natural order. In this case, for any given $\lambda_1 \geq 0, \lambda_2 \geq 0$, the function $\varphi$ takes the following form,

$$\varphi(y) = \lambda_1 \|y\|_1 + \lambda_2 \sum_{i=2}^{n} |y_i - y_{i-1}|. \tag{5.4}$$

The definition (5.4) contains the Lasso regularizer as a special case if we take $\lambda_2 = 0$.

By introducing an auxiliary variable $z \in \Re^n$, we can reformulate the fused Lasso problem into the framework of (1.1), i.e.,

$$\begin{aligned} \min_{y,y_0,z} \quad & f(y, y_0) + \varphi(z) \\ \text{s.t.} \quad & y - z = 0. \end{aligned} \tag{5.5}$$

The KKT system can be written as follows,

$$\nabla f(y, y_0) + [x; 0] = 0, \ z - \Pr_\varphi(x + z) = 0, \ y - z = 0.$$

The function $\widehat{\mathcal{L}}_\sigma(\cdot)$ in (1.7) can be specifically written as

$$\widehat{\mathcal{L}}_\sigma(y, y_0, z, x; \tilde{y}') = \hat{f}(\tilde{y}, \tilde{y}') + \varphi(z) + \langle y - z, x \rangle + \frac{\sigma}{2}\|y - z\|^2.$$

Consequently, the majorized iPADMM scheme for solving (5.5) can be described as follows,

$$\begin{cases} (y^{k+1}, y_0^{k+1}) = \arg\min_{y,y_0} \widehat{\mathcal{L}}_\sigma(y, y_0, z^k, x^k; \tilde{y}^k) + \frac{1}{2}\|\tilde{y} - \tilde{y}^k\|_{\mathcal{S}}^2, \\ z^{k+1} = \arg\min_z \varphi(z) + \frac{\sigma}{2}\|z - (y^{k+1} + x^k/\sigma)\|^2, \\ x^{k+1} = x^k + \tau\sigma(y^{k+1} - z^{k+1}), \end{cases} \tag{5.6}$$

where $\tau \in (0, (1+\sqrt{5})/2)$ and $\mathcal{S} = -\frac{1}{2}\widehat{\Sigma}_f + \text{Diag}(0, \sigma r)$ with $r$ being a given positive number. It is obviously that

$$\frac{1}{2}\widehat{\Sigma}_f + \mathcal{S} + \sigma \begin{pmatrix} I & 0 \\ 0 & 0 \end{pmatrix} \succ 0.$$

Based on the optimality conditions, we measure the accuracy of an approximate KKT point $(y, y_0, z, x)$ via

$$\eta_{\text{FL}} = \max\{\eta_P, \eta_D, \eta_C\}, \tag{5.7}$$

where

$$\eta_P = \left\{\frac{\|y - z\|}{1 + \|y\| + \|z\|}\right\}, \ \eta_D = \frac{\|\nabla f(y, y_0) + [x; 0]\|}{1 + \|\nabla f(y, y_0)\| + \|x\|}, \ \eta_C := \frac{\|z - \text{Pr}_\varphi(x + z)\|}{1 + \|x\| + \|z\|}.$$

It is worth mentioning that, the solution $z^{k+1}$ can be obtained by the following result.

**Lemma 5.1.** *For any $\lambda_1, \lambda_2 \geq 0$, the optimal solution $z^*$ of*

$$\min_z \lambda_1 \|z\|_1 + \lambda_2 \sum_{i=2}^n |z_i - z_{i-1}| + \frac{1}{2}\|z - v\|^2,$$

*can be described as*

$$z^* = \text{sgn}(z_0) \cdot \max(|z_0| - \lambda_1, 0),$$

*where $z_0 := \arg\min\limits_z \lambda_2 \sum_{i=2}^n |z_i - z_{i-1}| + \frac{1}{2}\|z - v\|^2$.*

The above result was first shown in [46] by using the subgradient technique and an alternative proof can be found in [39]. Though there is no closed-form expression of $z_0$ when $\lambda_2 > 0$, the algorithms presented in [6] can used to get $z_0$ efficiently. The corresponding code can download from website "*https://www.gipsa-lab.grenoble-inp.fr/ laurent.condat/software.html*".

By using [36, Theorem 1 & Remark 1] and [11, Theorem 3.1], we know that the KKT mapping $\mathcal{R}$ corresponding to (5.5) is metrically subregular at any KKT point for original when $\lambda_2 = 0$. Therefore, directly from Theorem 3.1, we can get the following results.

**Proposition 5.1.** *Let $\lambda_2 = 0$ and $\{u^k := (y^k, y_0^k, z^k, x^k)\}$ be the infinite sequence generated by the majorized iPADMM scheme (5.6). Then the sequence $\{u^k\}$ converges linearly to a KKT point of (5.5).*

## 5.2 Constrained logistic regression

Inspired by the work of James et al. [31], we consider the function $\varphi$ in the following form,

$$\varphi(y) = \lambda\|y\|_1 + \mathbb{I}_\mathcal{D}(y), \tag{5.8}$$

where $\mathcal{D} := \{y|\, Dy \geq d\}$, $D \in \Re^{m \times n}$, $d \in \Re^m$, function $\mathbb{I}_\mathcal{D}(\cdot)$ is an indicator function of convex set $\mathcal{D}$. We can rewrite (5.1) as

$$\begin{aligned} \min_{y, y_0, u, v} \quad & f(y, y_0) + \lambda\|v\|_1 + \mathbb{I}_{\Re^m_+}(u) \\ \text{s.t.} \quad & Dy - u = d, \\ & y - v = 0. \end{aligned} \tag{5.9}$$

The KKT conditions are given by

$$\nabla f(y, y_0) + [D^T\xi + \zeta; 0] = 0, \ v - \text{Pr}_{\lambda\|\cdot\|_1}(\zeta + v) = 0, \ u - \Pi_{\Re^m_+}(\xi + u) = 0, \ Du - u = d, \ y - v = 0,$$

and the majorized augmented Lagrangian function $\widehat{\mathcal{L}}_\sigma(\cdot)$ in (1.7) can be specifically written as

$$\widehat{\mathcal{L}}_\sigma(y, y_0, u, v, \xi, \zeta; \tilde{y}') = \hat{f}(\tilde{y}, \tilde{y}') + \lambda\|v\|_1 + \mathbb{I}_{\Re^m_+}(u) + \langle Dy - u - d, \xi \rangle + \langle y - v, \zeta \rangle$$
$$+ \frac{\sigma}{2}\|Dy - u - d\|^2 + \frac{\sigma}{2}\|y - v\|^2.$$

Therefore, we can solve (5.9) via the following iterative scheme:

$$\begin{cases} (y^{k+1}, y_0^{k+1}) = \arg\min_{y, y_0} \widehat{\mathcal{L}}_\sigma(y, y_0, u^k, v^k, \xi^k, \zeta^k; \tilde{y}^k) + \frac{1}{2}\|\tilde{y} - \tilde{y}^k\|_{\mathcal{S}}^2, \\ u^{k+1} = \max\{Dy^{k+1} - d + \xi^k/\sigma, \mathbf{0}\}, \\ v^{k+1} = \arg\min_v \lambda\|v\|_1 + \frac{\sigma}{2}\|v - (y^{k+1} + \zeta^k/\sigma)\|^2, \\ \xi^{k+1} = \xi^k + \tau\sigma(Dy^{k+1} - u^{k+1} - d), \\ \zeta^{k+1} = \zeta^k + \tau\sigma(y^{k+1} - v^{k+1}), \end{cases} \quad (5.10)$$

where $\tau \in (0, (1+\sqrt{5})/2)$ and $\mathcal{S} = -\frac{1}{2}\widehat{\Sigma}_f + \mathrm{Diag}(0, \sigma r)$ with $r$ being a given positive number. It is obviously that

$$\frac{1}{2}\widehat{\Sigma}_f + \mathcal{S} + \sigma \begin{pmatrix} D^T D + I & 0 \\ 0 & 0 \end{pmatrix} \succ 0.$$

Based on the optimality conditions, we measure the accuracy of an approximate KKT point $(y, y_0, u, v, \xi, \zeta)$ via

$$\eta_{\mathrm{CL}} = \max\{\eta_P, \eta_D, \eta_C\}, \quad (5.11)$$

where

$$\eta_P = \max\left\{\frac{\|Dy - u - d\|}{1 + \|Dy\| + \|u\| + \|d\|}, \frac{\|y - v\|}{1 + \|y\| + \|v\|}\right\}, \eta_D = \frac{\|\nabla f(y, y_0) + [D^T \xi + \zeta; 0]\|}{1 + \|\nabla f(y, y_0)\| + \|D^T \xi\| + \|\zeta\|},$$
$$\eta_C := \max\left\{\frac{\|v - \mathrm{Pr}_{\lambda\|\cdot\|_1}(\zeta + v)\|}{1 + \|\zeta\| + \|v\|}, \frac{\|u - \Pi_{\Re^m_+}(\xi + u)\|}{1 + \|\xi\| + \|u\|}.\right\}$$

## 5.3 Numerical experiments

In this subsection, we evaluate the performance of majorized iPADMM for solving Lasso logistic regression, fused Lasso logistic regression(5.5), constrained Lasso logistic regression (5.9), respectively. For all numerical experiments of ADMM-type method, we choose $\tau = 1.618$. All computational results are obtained by running Matlab R2015b on Mac OS X 10.10.5 (2.9 GHz Intel Core i5 16GB 1867 MHz DDR3).

Consider the following two self-adjoint linear operators:

$$\mathcal{S}_0 = \sigma \begin{pmatrix} 0 & 0 \\ 0 & r \end{pmatrix} \& \mathcal{S} = -\frac{1}{2}\widehat{\Sigma}_f + \sigma \begin{pmatrix} 0 & 0 \\ 0 & r \end{pmatrix}.$$

Note that the proximal term $\mathcal{S}$ may not be a positive semidefinite operator. In the subsequent discussions, we call the majorized iPADMM scheme with $\mathcal{S}_0$ as "majorized sPADMM". In all tests, we set $r = 10^{-6}$ and choose $\widehat{\Sigma}_f$ as (5.3) unless otherwise specified.

### 5.3.1 Lasso logistic regression: LIBSVM dataset

In this part, we apply the majorized iPADMM to Lasso logistic regression and test its performance with data sets: **a8a**, **a9a**, **colon-cancer**, **duke breast-cancer** (duke-BC). These data sets are obtained from the LIBSVM datasets [2].

Notice that the subproblem corresponding to $(y, y_0)$ can be reformulated as a linear system of equations. In order to exactly and efficiently solve this subproblem, we use different techniques to get the inverse matrix of $H := \widehat{\Sigma}_f + \sigma \text{Diag}(I, 0) + \mathcal{S}$. For **a8a, a9a**, i.e. $n \ll N$, $n$ is moderate, the inverse matrix $H \in \Re^{(n+1) \times (n+1)}$ can be efficiently obtained by MATLAB toolbox. For **conlon-cancer, duke-BC**, i.e. $n \gg N$, $N$ is moderate, the Sherman-Morrison-Woodbury formula [26] can be used to get the inverse matrix of $H$.

Figure 1 shows that the sequence generated by majorized iPADMM converges to a KKT point approximately as a linear rate. This is consistent with our Theorem 3.1. Moreover, we can also observe that the indefinite proximal term can improve the numerical performance. The majorized iPADMM brings about 35%– 60% reduction in the number of iterations needed as compared with the majorized sPADMM, when one chooses $\|u^k - \bar{u}\|_\mathcal{M} \leq 10^{-6}$.

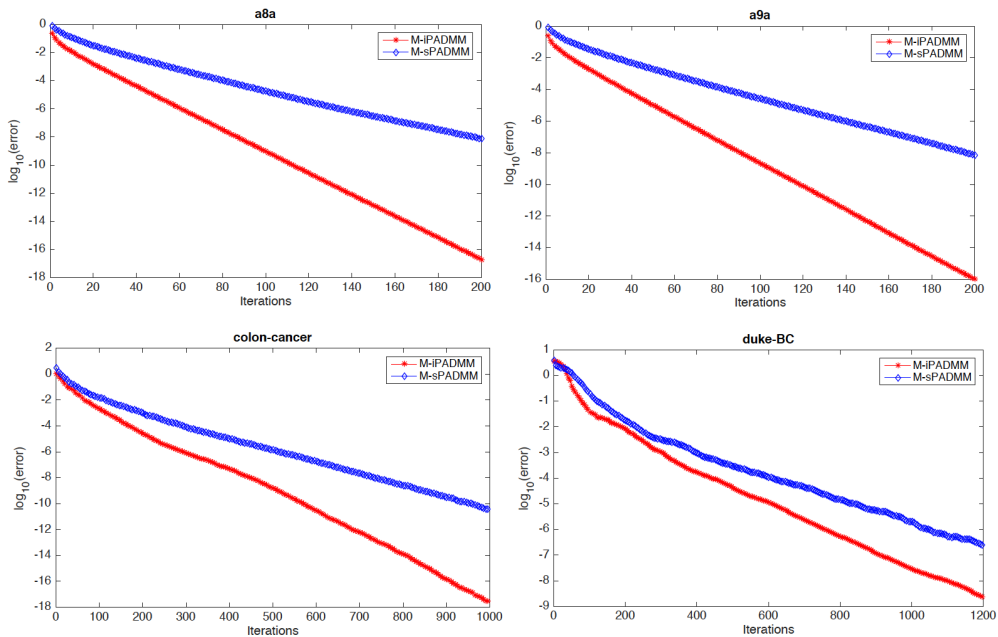

Figure 1: Comparison between the performance of our majorized iPADMM (M-iPADMM) and the majorized sPADMM (M-sPADMM) on datasets a8a, a9a, colon-cancer, and duke-BC. error := $\|u^k - \bar{u}\|_\mathcal{M}$.

### 5.3.2 Fused lasso logistic regression: LIBSVM dataset

In this part, we apply the majorized iPADMM to fused Lasso logistic regression and test its performance with data sets: **a8a**, **a9a**, **colon-cancer**, **duke breast-cancer**, **rcv1.binary**, **news20.binary**. These data sets are from the LIBSVM datasets [2].

In order to evaluate the numerical performance, we also report majorized sPADMM and a commonly used first order algorithm for solving fused lasso problem: the accelerated proximal gradient algorithm implemented in [39]. This algorithm is called "EFLA"[1] in [39]. In the comparison, we terminate "EFLA" if the objective value $\text{obj}_E$ obtained by "EFLA" satisfies

$$\frac{\text{obj}_E - \text{obj}_M}{|\text{obj}_M|} \leq 10\varepsilon,$$

where $\text{obj}_M$ is the objective value obtained by majorized iPADMM. To keep things fair, we use the algorithm in [6] to implement the subproblem in EFLA. This makes EFLA faster than the original code in software "SLEP 4.0".

In this test, we choose the regularization parameters as follows

$$\lambda_1 = \lambda_2 = \frac{\gamma}{N}\|B^T b\|_\infty,$$

where $N$ is the size of sample and $0 < \gamma < 1$. The majorized iPADMM and majorized sPADMM will be terminated when $\eta_{FL} < 10^{-6}$ or the maximum iteration number 50,000 is reached.

For **rcv1.binary, news20.binary**, i.e., both $n$ and $N$ are large, by taking $\widehat{\Sigma}_f := \sum_{i=1}^K \mu_i P_i P_i^T$, where $\sqrt{\mu_1}, \ldots, \sqrt{\mu_K}$ are the $K$ largest singular values of $A$ and $P_1, \ldots, P_K$ are the corresponding left-eigenvectors, the Sherman-Morrison-Woodbury formula can also be used to solve the linear system of equations corresponding to $(y, y_0)$ subproblem.

Now we are ready to report the comparison results. From Table 1, we can observe that the numerical performance of majorized iPADMM outperforms the other two methods for all cases. Besides, we can also see that the majorized iPADMM brings about 40%– 50% reduction in the number of iterations needed as compared with the majorized sPADMM except cases that **a8a** and **a9a** with $\gamma = 10^{-2}$.

| probName N \| n | $L_C$ | $\gamma$ | nnz | IterNum MiPA \| MsPA \| E | Time(sec) MiPA \| MsPA \| E |
|---|---|---|---|---|---|
| a8a 22696, 123 | 1.4e+05 | 1.0e-02 | 13 | 37 \| 48 \|106 | 0.10 \| 0.11 \| 0.31 |
| | | 1.0e-03 | 43 | 82 \| 168 \|331 | 0.18 \| 0.38 \| 0.88 |
| a9a 32561, 123 | 2.0e+05 | 1.0e-02 | 12 | 37 \| 46 \|100 | 0.12 \| 0.15 \| 0.42 |
| | | 1.0e-03 | 41 | 83 \| 171 \|291 | 0.24 \| 0.47 \| 1.13 |
| colon-cancer 62, 2000 | 1.9e+04 | 1.0e-02 | 52 | 319 \| 542 \|1294 | 0.34 \| 0.46 \| 0.47 |
| | | 1.0e-03 | 69 | 2010 \| 3891 \|6576 | 1.59 \| 3.05 \| 2.22 |
| duke-BC 44, 7129 | 1.1e+05 | 1.0e-02 | 61 | 750 \| 1271 \|3177 | 1.61 \| 2.72 \| 3.55 |
| | | 1.0e-03 | 68 | 2287 \| 6683 \|16105 | 4.84 \| 14.10 \| 18.17 |
| rcv1.binary 20242, 47236 | 4.5e+02 | 1.0e-02 | 123 | 650 \| 1099 \|1703 | 6.50 \| 11.15 \| 14.34 |
| | | 1.0e-03 | 744 | 2533 \| 4835 \|4640 | 25.18 \| 49.71 \| 38.96 |
| news20.binary 19996, 1355191 | 1.2e+03 | 1.0e-02 | 264 | 795 \| 1287 \|2618 | 178.52 \| 282.25 \| 278.54 |
| | | 1.0e-03 | 2973 | 2382 \| 4362 \|8860 | 536.45 \| 970.01 \| 935.00 |

Table 1: Comparision between the performance of majorized iPADMM (MiPA), majorized sPADMM (MsPA), and EPLA (E). "nnz" denotes the number of nonzeros in the solution $z$ generated by majorized iPADMM . "IterNum" denotes the number of iterations. $L_C := \lambda_{\max}(BB^T)$.

---

[1]The corresponding software SLEP 4.0 can download from website "http://yelab.net/software/SLEP/"

### 5.3.3 Constrained logistic regression: synthetic data

This part tests the performance of majorized iPADMM for the constrained logistic regression by using synthetic data. The data $B \in \Re^{n \times N}$, $D \in \Re^{m \times n}$ and $d \in \Re^m$ are generated by the standard normal distribution. This example can also be used to illustrate the claim mentioned in the introduction that our majorized iPADMM with the positive semidefinite operator $\widehat{\Sigma}_f$ outperforms that with the Lipschitz constant $\lambda_{\max}(\widehat{\Sigma}_f)$ (denoted as L-majorized ADMM).

We compare our majorized iPADMM with majorized sPADMM, L-majorized ADMM with an indefinite proximal term $(\mathrm{Diag}(-\frac{1}{2}\lambda_{\max}(\widehat{\Sigma}_f)I, \sigma r))$, and L-majorized sPADMM. All the algorithms will be terminated when $\eta_{CL} < 10^{-5}$ or the maximum iteration number 50,000 is reached. In this test, we choose the regularization parameter $\lambda = \gamma \|B^T b\|_\infty / N$, where $0 < \gamma < 1$.

Table 2 reports the number of iterations, runtime of four different methods. From the table, we can see that our majorized iPADMM outperforms all the other three methods. In each case, the majorized iPADMM can sometimes bring about 40% reduction in the number of iterations needed for convergence as compared with the majorized sPADMM. Note that though the size of each scenario is small, some cases still can not be solved within maximum iteration by using L-majorized iPADMM and L-majorized sPADMM.

Figure 2 shows the relative KKT residual norm $\eta_{CL}$ of the iterates generated by majorized iPADMM and L-majorized iPADMM on two different synthetic data sets. This figure can be used to illustrate the advantage of using the positive semidefinite operator $\widehat{\Sigma}_f$ over the Lipschitz constant $\lambda_{\max}(\widehat{\Sigma}_f)$.

| N,n,m | $L_C$ | $\lambda$ | IterNum<br>MiPA \| MsPA \| LiA \| LA | Time (sec)<br>MiPA \| MsPA \| LiA \| LA |
|---|---|---|---|---|
| 2000,500,50 | 2.25 | $10^{-2}$ | 306.7 \| 307.3 \| 307.5 \| 312.7 | 1.5 \| 1.5 \| 1.5 \| 1.5 |
| | | $10^{-3}$ | 154.3 \| 161.3 \| 167.3 \| 186.4 | 0.8 \| 0.9 \| 0.9 \| 1.0 |
| | | $10^{-4}$ | 115.2 \| 119.8 \| 128.1 \| 149.9 | 0.6 \| 0.6 \| 0.7 \| 0.8 |
| 30,50,20 | 4.59 | $10^{-2}$ | 317.0 \| 374.0 \| 1917.3 \| 3825.7 | 0.1 \| 0.1 \| 0.3 \| 0.5 |
| | | $10^{-3}$ | 2014.9 \| 3887.8 \| 18494.4 \| 35075.7 | 0.2 \| 0.4 \| 1.9 \| 3.5 |
| | | $10^{-4}$ | 11636.8 \| 21473.3 \| 49961.2 \| 50000.0 | 1.0 \| 1.6 \| 4.9 \| 4.7 |
| 50,200,30 | 9.08 | $10^{-2}$ | 573.3 \| 581.8 \| 3470.5 \| 6734.6 | 0.4 \| 0.4 \| 1.8 \| 3.3 |
| | | $10^{-3}$ | 1550.8 \| 2303.0 \| 27501.7 \| 47114.9 | 0.8 \| 1.2 \| 12.8 \| 21.8 |
| | | $10^{-4}$ | 6403.4 \| 12001.4 \| 50000.0 \| 50000.0 | 3.2 \| 5.7 \| 23.2 \| 23.2 |

Table 2: Comparision between the performance of majorized iPADMM (MiPA), majorized sPADMM (MsPA), L-majorized iPADMM (LiA), and L-majoirized sPADMM (LA); "IterNum" denotes the number of iterations. $L_C := \lambda_{\max}(AA^T)/N$. All results are averaged over 10 instances.

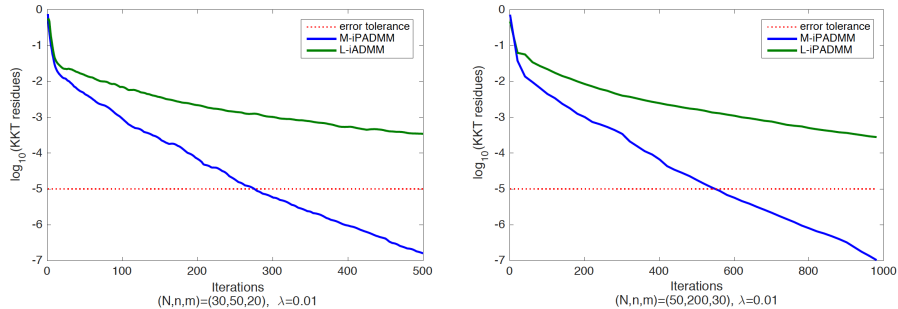

Figure 2: Comparision between the performance of majorized iPADMM (M-iPADMM) and L-majorized iPADMM (L-iPADMM) on synthetic data within given iteration numbers. All results are averaged over 10 instances.

# 6 Conclusion remarks

In this paper, we have established the linear rate convergence of the majorized ADMM with indefinite proximal terms for solving the 2-block linearly constrained convex composite optimization problem under a metric subregularity assumption. Numerical results on three types of regularized logistic regression have been given to evaluate the effectiveness of the 2-block majorized ADMM with indefinite proximal terms. From these results, we can see that, for many cases, the majorized ADMM with indefinite proximal terms can bring about 30%–50% reduction in the number of iterations needed for convergence as compared with the majorized ADMM with semi-proximal terms.

Strongly motivated by the numerical performance of the symmetric Gauss-Seidel based ADMM for solving multi-block convex composite quadratic programming, we also proved the linear rate of convergence of a symmetric Gauss-Seidel based majorized ADMM with indefinite proximal terms by building its equivalence to the 2-block majorized ADMM with specially constructed proximal terms (possibly indefinite). This will extremely facilitate the future exploration of the multi-block general linear/nonlinear models. We leave this topic as our future work.

# A   Appendix: Proofs of the Lemmas

## A.1   Proof of Lemma 3.1

*Proof.* The optimality condition for (1.8a) is

$$0 \in \partial p(y^{k+1}) + \nabla f(y^k) + \mathcal{A}x^k + \sigma\mathcal{A}(\mathcal{A}^*y^{k+1} + \mathcal{B}^*z^k - c) + (\mathcal{S} + \widehat{\Sigma}_f)(y^{k+1} - y^k). \tag{A.1}$$

It follows from (1.8c) that

$$x^k + \sigma(\mathcal{A}^*y^{k+1} + \mathcal{B}^*z^k - c) = \tau^{-1}(x^{k+1} - x^k) + x^k - \sigma\mathcal{B}^*(z^{k+1} - z^k). \tag{A.2}$$

From Proposition 2.1, we know that there exists $W_f^k \in \text{conv}\partial^2 f[y^{k+1}, y^k]$ such that

$$\nabla f(y^k) - \nabla f(y^{k+1}) = W_f^k(y^k - y^{k+1}).$$

Substituting the above equation and (A.2) into (A.1), we get

$$0 \in \partial p(y^{k+1}) + \nabla f(y^{k+1}) + \mathcal{A}[x^k + \tau^{-1}(x^{k+1} - x^k) - \sigma\mathcal{B}^*(z^{k+1} - z^k)] + \mathcal{S}^k(y^{k+1} - y^k),$$

where $\mathcal{S}^k := \mathcal{S} + \widehat{\Sigma}_f - W_f^k$. Thus, it holds that

$$y^{k+1} = \text{Pr}_p\{y^{k+1} - \nabla f(y^{k+1}) - \mathcal{A}[x^k + \tau^{-1}(x^{k+1} - x^k) - \sigma\mathcal{B}^*(z^{k+1} - z^k)] - \mathcal{S}^k(y^{k+1} - y^k)\}.$$

Similarly, there exists $W_g^k \in \text{conv}\partial^2 g(z^{k+1}, z^k)$ such that $\nabla g(z^k) - \nabla g(z^{k+1}) = W_g^k(z^k - z^{k+1})$ and

$$z^{k+1} = \text{Pr}_q\{z^{k+1} - \nabla g(z^{k+1}) - \mathcal{B}[\tau^{-1}(x^{k+1} - x^k) + x^k] - \mathcal{T}^k(z^{k+1} - z^k)\},$$

where $\mathcal{T}^k := \mathcal{T} + \widehat{\Sigma}_g - W_g^k$. Since the Moreau-Yosida proximal mappings $\text{Pr}_p(\cdot)$ and $\text{Pr}_q(\cdot)$ are globally Lipschitz continuous, one has that for any $k \geq 1$,

$$\begin{aligned}
&\|\mathcal{R}(u^{k+1})\|^2 \\
&\leq \|\mathcal{S}^k(y^{k+1} - y^k) - \sigma\mathcal{A}\mathcal{B}^*(z^{k+1} - z^k) + (\tau^{-1} - 1)\mathcal{A}(x^{k+1} - x^k)\|^2 \\
&\quad + \|\mathcal{T}^k(z^{k+1} - z^k) + (\tau^{-1} - 1)\mathcal{B}(x^{k+1} - x^k)\|^2 + \|(\tau\sigma)^{-1}(x^{k+1} - x^k)\|^2 \\
&\leq 3\lambda_{\max}^2(|\mathcal{S}^k|)\|y^{k+1} - y^k\|^2 + 3\sigma\lambda_{\max}(\mathcal{A}^*\mathcal{A})\|(z^{k+1} - z^k)\|_{\sigma\mathcal{B}\mathcal{B}^*}^2 \\
&\quad + 2\lambda_{\max}^2(|\mathcal{T}^k|)\|z^{k+1} - z^k\|^2 + \|(\tau\sigma)^{-1}(x^{k+1} - x^k)\|^2 \\
&\quad + (1 - \tau^{-1})^2[2\lambda_{\max}(\mathcal{B}^*\mathcal{B})\|(x^{k+1} - x^k)\|^2 + 3\lambda_{\max}(\mathcal{A}^*\mathcal{A})\|(x^{k+1} - x^k)\|^2].
\end{aligned} \tag{A.3}$$

Next, we estimate upper bounds of $\lambda_{\max}(\mathcal{S}^k)$ and $\lambda_{\max}(\mathcal{T}^k)$, respectively. It follows from $\Sigma_f \preceq W_f^k \preceq \widehat{\Sigma}_f$ and $\Sigma_g \preceq W_g^k \preceq \widehat{\Sigma}_g$ that

$$-\frac{1}{2}\widehat{\Sigma}_f \preceq \mathcal{S}^k \preceq \mathcal{S} + \widehat{\Sigma}_f \text{ and } -\frac{1}{2}\widehat{\Sigma}_g \preceq \mathcal{T}^k \preceq \mathcal{T} + \widehat{\Sigma}_g, \ \forall k \geq 1.$$

Then

$$\lambda_{\min}(\mathcal{S}^k) \geq -\tfrac{1}{2}\lambda_{\max}(\widehat{\Sigma}_f) \text{ and } \lambda_{\max}(\mathcal{S}^k) \leq \lambda_{\max}(\mathcal{S} + \tfrac{1}{2}\widehat{\Sigma}_f) + \tfrac{1}{2}\lambda_{\max}(\widehat{\Sigma}_f),$$

$$\lambda_{\min}(\mathcal{T}^k) \geq -\tfrac{1}{2}\lambda_{\max}(\widehat{\Sigma}_g) \text{ and } \lambda_{\max}(\mathcal{T}^k) \leq \lambda_{\max}(\mathcal{T} + \tfrac{1}{2}\widehat{\Sigma}_g) + \tfrac{1}{2}\lambda_{\max}(\widehat{\Sigma}_g),$$

and consequently, for any $k \geq 1$, one has

$$\lambda_{\max}(\mathcal{S}^k) \leq \lambda_{\max}(\mathcal{S} + \tfrac{1}{2}\widehat{\Sigma}_f) + \tfrac{1}{2}\lambda_{\max}(\widehat{\Sigma}_f), \lambda_{\max}(\mathcal{T}^k) \leq \lambda_{\max}(\mathcal{T} + \tfrac{1}{2}\widehat{\Sigma}_g) + \tfrac{1}{2}\lambda_{\max}(\widehat{\Sigma}_g).$$

By substituting the above two inequalities into (A.3), we can get (3.4). This completes the proof. □

## A.2 Proof of Lemma 3.2

*Proof.* For any $\varepsilon \in \Re$, define
$$h_\varepsilon(x) := h(x) + \frac{\varepsilon^2}{2}\|x\|^2, \ \forall\, x \in \mathcal{X}.$$

Then, similar to the proof of [41, Theorem 2.1.5], one has
$$\langle \nabla h_\varepsilon(x) - \nabla h_\varepsilon(\bar{x}), x - \bar{x}\rangle \geq \|\nabla h_\varepsilon(x) - \nabla h_\varepsilon(\bar{x})\|^2_{(\mathcal{P}+\varepsilon^2\mathcal{I})^{-1}}, \ \forall\, \varepsilon \neq 0.$$

Consequently, for any $\varepsilon \neq 0$, it holds that
$$\begin{aligned}
&\langle \nabla h_\varepsilon(x) - \nabla h_\varepsilon(\bar{x}), y - \bar{x}\rangle \\
&\geq \|\nabla h_\varepsilon(x) - \nabla h_\varepsilon(\bar{x})\|^2_{(\mathcal{P}+\varepsilon^2\mathcal{I})^{-1}} + \langle \nabla h_\varepsilon(x) - \nabla h_\varepsilon(\bar{x}), y - x\rangle \\
&= \|(\mathcal{P}+\varepsilon^2\mathcal{I})^{-1/2}(\nabla h_\varepsilon(x) - \nabla h_\varepsilon(\bar{x})) + \tfrac{1}{2}(\mathcal{P}+\varepsilon^2 I)^{1/2}(y-x)\|^2 - \tfrac{1}{4}\|x-y\|^2_{(\mathcal{P}+\varepsilon^2\mathcal{I})} \\
&\geq -\tfrac{1}{4}\|x-y\|^2_{(\mathcal{P}+\varepsilon^2\mathcal{I})}.
\end{aligned}$$

This together with the definition of $h$ implies that
$$\langle \nabla h(x) - \nabla h(\bar{x}), y - \bar{x}\rangle + \langle \varepsilon^2(x - \bar{x}), y - \bar{x}\rangle \geq -\frac{1}{4}\|x-y\|^2_{(\mathcal{P}+\varepsilon^2\mathcal{I})}, \ \forall \varepsilon \neq 0.$$

Therefore, by taking limits on both sides of the above inequality for $\varepsilon \to 0$, we complete the proof. $\square$

## A.3 Proof of Lemma 3.3

*Proof.* By the first order optimality conditions of (1.8a) and (1.8b), one has
$$\begin{cases} 0 \in \partial p(y^{k+1}) + \nabla f(y^k) + (\widehat{\Sigma}_f + \mathcal{S})(y^{k+1} - y^k) + \mathcal{A}(x^k + \sigma(\mathcal{A}^*y^{k+1} + \mathcal{B}^*z^k - c)), \\ 0 \in \partial q(z^{k+1}) + \nabla g(z^k) + (\widehat{\Sigma}_g + \mathcal{T})(z^{k+1} - z^k) + \mathcal{B}(x^k + \sigma r^{k+1}). \end{cases} \quad (A.4)$$

Since $(\bar{y}, \bar{z}, \bar{x})$ is a KKT point, it holds that
$$\begin{cases} 0 \in \partial p(\bar{y}) + \nabla f(\bar{y}) + \mathcal{A}\bar{x}, \\ 0 \in \partial q(\bar{z}) + \nabla g(\bar{z}) + \mathcal{B}\bar{x}. \end{cases} \quad (A.5)$$

It follows from the maximal monotonicity of $\partial p$ that
$$\begin{aligned}
0 \leq &\langle -\mathcal{A}(x^k + \sigma(\mathcal{A}^*y^{k+1} + \mathcal{B}^*z^k - c)) + \mathcal{A}\bar{x}, y^{k+1} - \bar{y}\rangle + \langle \nabla f(\bar{y}) - \nabla f(y^k), y^{k+1} - \bar{y}\rangle \\
&- \langle (\widehat{\Sigma}_f + \mathcal{S})(y^{k+1} - y^k), y^{k+1} - \bar{y}\rangle.
\end{aligned}$$

Thus, by reorganizing the above inequality and Lemma 3.2, one has
$$\begin{aligned}
&\langle \bar{x} - (x^k + \sigma(\mathcal{A}^*y^{k+1} + \mathcal{B}^*z^k - c)), \mathcal{A}^*(y^{k+1} - \bar{y})\rangle - \langle (\widehat{\Sigma}_f + \mathcal{S})(y^{k+1} - y^k), y^{k+1} - \bar{y}\rangle \\
&\geq \langle \nabla f(y^k) - \nabla f(\bar{y}), y^{k+1} - \bar{y}\rangle \geq -\tfrac{1}{4}\|y^{k+1} - y^k\|^2_{\widehat{\Sigma}_f}.
\end{aligned} \quad (A.6)$$

Similarly, by using the maximal monotonicity of $\partial q$ and Lemma 3.2, it holds that

$$\langle \bar{x} - (x^k + \sigma r^{k+1}), \mathcal{B}^*(z^{k+1} - \bar{z})\rangle - \langle (\widehat{\Sigma}_g + \mathcal{T})(z^{k+1} - z^k), z^{k+1} - \bar{z}\rangle \\ \geq -\tfrac{1}{4}\|z^{k+1} - z^k\|^2_{\widehat{\Sigma}_g}. \tag{A.7}$$

By adding (A.6) and (A.7) together, we have

$$\Delta_k + \langle(\widehat{\Sigma}_f + \mathcal{S})(y^{k+1} - y^k), \bar{y} - y^{k+1}\rangle + \langle(\widehat{\Sigma}_g + \mathcal{T})(z^{k+1} - z^k), \bar{z} - z^{k+1}\rangle \\ \geq -\tfrac{1}{4}\|y^{k+1} - y^k\|^2_{\widehat{\Sigma}_f} - \tfrac{1}{4}\|z^{k+1} - z^k\|^2_{\widehat{\Sigma}_g}, \tag{A.8}$$

where

$$\Delta_k := \langle \bar{x} - (x^k + \sigma r^{k+1}) + \sigma \mathcal{B}^*(z^{k+1} - z^k), \mathcal{A}^*(y^{k+1} - \bar{y})\rangle + \langle \bar{x} - (x^k + \sigma r^{k+1}), \mathcal{B}^*(z^{k+1} - \bar{z})\rangle.$$

Directly from [34, Equation (47)], it holds that

$$\begin{aligned}\Delta_k = \ &(2\tau\sigma)^{-1}\left(\|x^k - \bar{x}\|^2 - \|x^{k+1} - \bar{x}\|^2\right) + (2\sigma)^{-1}(\tau - 1)\|r^{k+1}\|^2 \\ &- \tfrac{\sigma}{2}\|\mathcal{A}^*y^{k+1} + \mathcal{B}^*z^k - c\|^2 + \tfrac{\sigma}{2}\left(\|\mathcal{B}^*z^k - \mathcal{B}^*\bar{z}\|^2 - \|\mathcal{B}^*z^{k+1} - \mathcal{B}^*\bar{z}\|^2\right).\end{aligned} \tag{A.9}$$

Since for any self-adjoint linear operator $\mathcal{G}$, it holds that $\langle u, \mathcal{G}v\rangle = \tfrac{1}{2}(\|u+v\|^2_\mathcal{G} - \|u\|^2_\mathcal{G} - \|v\|^2_\mathcal{G})$, we have

$$\begin{aligned}&\langle(\widehat{\Sigma}_f + \mathcal{S})(y^{k+1} - y^k), \bar{y} - y^{k+1}\rangle + \langle(\widehat{\Sigma}_g + \mathcal{T})(z^{k+1} - z^k), \bar{z} - z^{k+1}\rangle \\ &= \tfrac{1}{2}\left(\|y^k - \bar{y}\|^2_{\widehat{\Sigma}_f + \mathcal{S}} - \|y^{k+1} - \bar{y}\|^2_{\widehat{\Sigma}_f + \mathcal{S}}\right) - \tfrac{1}{2}\|y^{k+1} - y^k\|^2_{\widehat{\Sigma}_f + \mathcal{S}} \\ &\quad + \tfrac{1}{2}\left(\|z^k - \bar{z}\|^2_{\widehat{\Sigma}_g + \mathcal{T}} - \|z^{k+1} - \bar{z}\|^2_{\widehat{\Sigma}_g + \mathcal{T}}\right) - \tfrac{1}{2}\|z^{k+1} - z^k\|^2_{\widehat{\Sigma}_g + \mathcal{T}}.\end{aligned}$$

This, together with (A.9) and (A.8) implies that the conclusion holds. The proof is completed. □

### A.4 Proof of Lemma 3.4

*Proof.* From [34, Lemma 7], it holds that

$$\begin{aligned}&\|y^{k+1} - y^k\|^2_{\frac{1}{2}\widehat{\Sigma}_f + \mathcal{S}} + \|z^{k+1} - z^k\|^2_{\frac{1}{2}\widehat{\Sigma}_g + \mathcal{T}} + \sigma\|\mathcal{A}^*y^{k+1} + \mathcal{B}^*z^k - c\|^2 + (1-\tau)\sigma\|r^{k+1}\|^2 \\ &\geq \|y^{k+1} - y^k\|^2_{\frac{1}{2}\widehat{\Sigma}_f + \mathcal{S}} + \|z^{k+1} - z^k\|^2_{\frac{1}{2}\widehat{\Sigma}_g + \mathcal{T}} + \|z^{k+1} - z^k\|^2_{\widehat{\Sigma}_g + \mathcal{T}} - \|z^k - z^{k-1}\|^2_{\widehat{\Sigma}_g + \mathcal{T}} \\ &\quad + \min(\tau, 1 + \tau - \tau^2)\sigma(\tau^{-1}\|r^{k+1}\|^2 + \|\mathcal{B}^*(z^{k+1} - z^k)\|^2) + (1 - \min(\tau, \tau^{-1}))\sigma(\|r^{k+1}\|^2 - \|r^k\|).\end{aligned}$$

This, together with the definition of $\phi_k$ and Lemma 3.3 implies the conclusion. □

### A.5 Proof of Lemma 3.5

*Proof.* First we show that

$$\tfrac{1}{2}\widehat{\Sigma}_f + \mathcal{S} + \sigma\mathcal{A}\mathcal{A}^* \succ 0 \ \& \ \tfrac{1}{2}\widehat{\Sigma}_g + \mathcal{T} + \sigma\mathcal{B}\mathcal{B}^* \succ 0 \Rightarrow \mathcal{H} \succ 0, \ \mathcal{M} \succ 0.$$

Suppose that $\frac{1}{2}\widehat{\Sigma}_f + \mathcal{S} + \sigma\mathcal{A}\mathcal{A}^* \succ 0$ & $\frac{1}{2}\widehat{\Sigma}_g + \mathcal{T} + \sigma\mathcal{B}\mathcal{B}^* \succ 0$, but there exists a vector $0 \neq d := (d_y, d_z, d_x) \in \mathcal{Y} \times \mathcal{Z} \times \mathcal{X}$ such that $\langle d, \mathcal{H}d \rangle = 0$, by using the definition of $\mathcal{H}$, we have

$$\langle d_y, \mathcal{H}_f d_y \rangle = 0, \ \langle d_z, \mathcal{H}_g d_z \rangle = 0, \ d_x = 0, \ \mathcal{E}^*(d_y, d_z, d_x) = 0.$$

Since $t_\tau > 0$ and $\frac{1}{2}\widehat{\Sigma}_g + \mathcal{T} + \sigma\mathcal{B}\mathcal{B}^* \succ 0$, we know $d_z = 0$. Consequently, $\mathcal{A}^* d_y = 0$. This together with the assumption that $\frac{1}{2}\widehat{\Sigma}_f + \mathcal{S} + \sigma\mathcal{A}\mathcal{A}^* \succ 0$ implies $d_y = 0$. This contradiction shows that $\mathcal{H} \succ 0$. We can get $\mathcal{M} \succ 0$ by using the same techniques. For brevity, we omit the details.

Next, we show that $\mathcal{H} \succ 0 \Rightarrow \frac{1}{2}\widehat{\Sigma}_f + \mathcal{S} + \sigma\mathcal{A}\mathcal{A}^* \succ 0$ & $\frac{1}{2}\widehat{\Sigma}_g + \mathcal{T} + \sigma\mathcal{B}\mathcal{B}^* \succ 0$. Since $t_\tau > 0$ and for any $d = (d_y, 0, 0) \in \mathcal{Y} \times \mathcal{Z} \times \mathcal{X}$, we have $\langle d, \mathcal{H}d \rangle = \langle d_y, (\mathcal{H}_f + \frac{1}{4}t_\tau \sigma \mathcal{A}\mathcal{A}^*)d_y \rangle$, then $\frac{1}{2}\widehat{\Sigma}_f + \mathcal{S} + \sigma\mathcal{A}\mathcal{A}^* \succ 0$ by the definition of $\mathcal{H}_f$. Similarly, for any $d = (0, d_z, 0) \in \mathcal{Y} \times \mathcal{Z} \times \mathcal{X}$, $\langle d, \mathcal{H}d \rangle = \langle d_z, (\mathcal{H}_g + \frac{1}{4}t_\tau \sigma \mathcal{B}\mathcal{B}^*)d_z \rangle$, then we can obtain that $\frac{1}{2}\widehat{\Sigma}_g + \mathcal{T} + \sigma\mathcal{B}\mathcal{B}^* \succ 0$. The proof is completed. □